\documentclass[amssymb,amsfonts,refcheck,11pt,verbatim,righttag]{amsart}
\usepackage{amsfonts}
\usepackage{amssymb}
\usepackage{amsthm}

 \numberwithin{equation}{section}
\theoremstyle{plain}

\numberwithin{equation}{section}
\newcommand{\Proof}{\noindent {\it Proof.~~}}

\newtheorem{thm}{Theorem}[section]

\newtheorem{lem}[thm]{Lemma}
\newtheorem{pro}[thm]{Proposition}
\newtheorem{cor}[thm]{Corollary}

\newtheorem{de}[thm]{Definition}
\newtheorem{rem}[thm]{Remark}

\newcommand{\R} {{\Bbb R}}
\newcommand{\N} {{\Bbb N}}
\newcommand{\Z} {{\Bbb Z}}
\newcommand{\M} {{\mathcal M}}

\newcommand{\I} {{\mathcal I}}
\newcommand{\A} {{\mathcal A}}

\newcommand{\C} {{\mathcal C}}
\newcommand{\D} {{\mathcal D}}

\newcommand{\ba}{{\bf a}}
\newcommand{\bb} {{\bf b}}

\newcommand{\lan}{{\mathcal L}}

\begin{document}
\baselineskip 16pt

\title{Weighted equilibrium states  for factor maps between subshifts}

\author{De-Jun Feng}
\address{
Department of Mathematics\\
The Chinese University of Hong Kong\\
Shatin,  Hong Kong\\
}
\email{djfeng@math.cuhk.edu.hk}
\date{}

\begin{abstract}
Let $\pi:X\to Y$ be a factor map, where $(X,\sigma_X)$ and  $(Y,\sigma_Y)$ are subshifts over finite alphabets.   Assume that
$X$ satisfies weak specification. Let $\ba=(a_1,a_2)\in \R^2$ with $a_1>0$ and $a_2\geq 0$. Let $f$ be a continuous function on $X$ with sufficient regularity (H\"{o}lder continuity, for instance). We show that there is a unique shift invariant  measure $\mu$ on $X$ that maximizes $\mu(f)+a_1h_\mu(\sigma_X)+ a_2h_{\mu\circ \pi^{-1}}(\sigma_Y)$. In particular, taking $f\equiv 0$ we see that there is a unique invariant measure $\mu$ on $X$ that maximizes the weighted entropy $a_1h_\mu(\sigma_X)+ a_2h_{\mu\circ \pi^{-1}}(\sigma_Y)$. This answers an open question raised by Gatzouras and Peres in \cite{GaPe96}. An extension is also given to  high dimensional cases. As an application, we show the  uniqueness of invariant measures with full Hausdorff dimension for certain  affine invariant sets on the $k$-torus under a diagonal endomorphism.

\end{abstract}

\maketitle

\section{Introduction}

%\end{document}

 Let $k\geq 2$ be an integer. Assume that $(X_i,\sigma_{X_i})$, $i=1,\ldots,k$,  are one-sided (or two-sided) subshifts over finite alphabets. Furthermore assume that $X_{i+1}$ is a factor of $X_i$ with a factor map $\pi_i:\; X_i\to X_{i+1}$ for $i=1,\ldots,k-1$.
 For convenience,  we use $\pi_0$ to denote the identity map on $X_1$.
  Define
 $\tau_i:\;X_1\to X_{i+1}$ by $\tau_i=\pi_i\circ\pi_{i-1}\circ
\cdots \circ \pi_0$ for $i=0,1,\ldots,k-1$.
Let $\M(X_i,\sigma_{X_i})$ denote the set of all $\sigma_{X_i}$-invariant Borel probability measures on $X_i$, endowed with the weak-star topology. For $f\in C(X_1)$ (the set of continuous functions on $X_1$),  and $\ba=(a_1,a_2,\ldots,a_k)\in \R^k$ with $a_1>0$ and $a_i\geq 0$ for $i\geq 2$, we say that $\mu\in \M(X_1,\sigma_{X_1})$ is an {\it $\ba$-weighted equilibrium state of $f$ for the factor maps $\pi_i$'s}, or simply, {\it $\ba$-weighted equilibrium state of $f$} if
\begin{equation}
\label{e-uu}
\mu(f) +\sum_{i=1}^ka_ih_{\mu\circ \tau_{i-1}^{-1}}(\sigma_{X_i})=\sup_{\eta\in \M(X,\sigma_X)}\left(\eta(f) +\sum_{i=1}^ka_ih_{\eta\circ \tau_{i-1}^{-1}}(\sigma_{X_i})\right),
\end{equation}
where $\mu(f)=\int_{X_1} f\; d\mu$, $\mu\circ\tau_{i-1}^{-1}$ denotes the measure on $X_{i}$ given by $\mu\circ \tau_{i-1}^{-1}(B)=\mu(\tau_{i-1}^{-1}(B))$ for any Borel set $B\subseteq X_{i}$, $h_{\mu\circ \tau_{i-1}^{-1}} (\sigma_{X_i})$  denotes the measure theoretic entropy of $\mu\circ \tau_{i-1}^{-1}$.
The right hand side of \eqref{e-uu} is called the {\it $\ba$-weighted topological pressure of $f$} and is denoted by $P^\ba(\sigma_{X_1},f)$.  The existence of at least one $\ba$-weighted equilibrium measure follows from the upper semi-continuity of the entropy functions $h_{(\cdot)}(\sigma_{X_i})$. In this paper we want to give conditions on $f$ and $X_i$'s to guarantee a unique $\ba$-weighted equilibrium state. The question  seems  quite fundamental in ergodic theory and symbolic
dynamics.

   We say that $X_1$ satisfies {\it weak specification} if there exists $p\in \N$ such that, for any two words $I$ and $J$ that are legal in $X_1$ (i.e., may be extended to sequences in $X_1$), there is a word $K$ of length not exceeding $ p$ such that the word $IKJ$ is legal in $X_1$. Similarly, say that $X_1$ satisfies {\it specification} if there exists $p\in \N$ such that, for any two words $I$ and $J$ that are legal in $X_1$, there is a word $K$ of length $p$ such that the word $IKJ$ is legal in $X_1$. For more details about the definitions, see Sect.~\ref{S-2}.

   For $f\in C(X_1)$ and $n\geq 1$ let
  \begin{equation}
  \label{e-sn}
  S_nf(x)=\sum_{i=0}^{n-1}f(\sigma^i_{X_1}x),\quad x\in X_1.
  \end{equation}
  Let $V(\sigma_{X_1})$ denote the set of $f\in C(X_1)$ such that there exists $c>0$ such that
  \begin{equation}
  \label{e-sn1}
  |S_nf(x)-S_nf(y)|\leq c\quad \mbox{ whenever } x_i=y_i \mbox{ for all } 0<i\leq n.
  \end{equation}
 Endow $X_1$ with the usual metric (see Sect.~\ref{S-2}). Clearly $V(\sigma_{X_1})$ contains  all H\"{o}lder continuous functions on $X_1$.
 The main result of the paper is the following.

\begin{thm}
 \label{thm-1}
  Assume that $X_1$ satisfies weak specification. Then for any $f\in V(\sigma_{X_1})$ and  $\ba=(a_1,a_2,\ldots, a_k)\in \R^k$ with   $a_1>0$ and $a_i\geq 0$ for $i\geq 2$,  $f$ has a unique $\ba$-weighted equilibrium state $\mu$. The measure $\mu$ is ergodic and, there exist $p\in \N$ and $c>0$ such that
$$\liminf_{n\to \infty}\sum_{i=0}^p\mu(A\cap \sigma^{-n-i}_{X_1}(B))\geq c\mu(A)\mu(B), \quad \forall\mbox{ Borel sets }A,B\subseteq X_1.
$$
Furthermore, if $X_1$ satisfies specification, then  there exists $c>0$ such that
  $$\liminf_{n\to \infty}\mu(A\cap \sigma^{-n}_{X_1}(B))\geq c\mu(A)\mu(B), \quad \forall\mbox{ Borel sets }A,B\subseteq X_1.
$$
  \end{thm}

When $\ba=(1,0,\ldots,0)$ and $X_1$ satisfies specification, Theorem \ref{thm-1} reduces to Bowen's theory   about the uniqueness of classical equilibrium states for the  subshift case \cite{Bow74}.  Taking $f=0$ in Theorem \ref{thm-1} yields,  whenever $X_1$ satisfies weak specification,  there
is a unique $\sigma_{X_1}$-invariant measure $\mu$ which maximizes the $\ba$-weighted entropy $\sum_{i=1}^ka_ih_{\mu\circ \tau_{i-1}^{-1}}(\sigma_{X_i})$.
Since
each irreducible  subshift of finite type satisfies weak specification (cf. Sect.~\ref{S-2}), this solves the following  open  question raised by Gatzouras and Peres (see
\cite[Problem 3]{GaPe96}):

{\it  Let  $\pi:\; X\to Y$ be a factor map between subshifts $X$ and $Y$, where  $X$ is an irreducible subshift of finite type. Let $\alpha>0$. Is there  a unique invariant measure $\mu$ maximizing the weighted entropy $h_\mu(\sigma_X)+\alpha h_{\mu\circ \pi^{-1}}(\sigma_Y)$? }

The above question  is closely related to dimension theory of  non-conformal dynamical systems. Let $T$ be the endmorphism on the $k$-dimensional torus ${\Bbb T}^k=\R^k/\Z^k$
 represented by an integral diagonal matrix $$\Lambda={\rm diag}(m_1,m_2,\ldots, m_k),$$ where $2\leq m_1\leq \ldots\leq  m_k.$
Let $\A$ denote the Cartesian product $$\prod_{i=1}^k\{0,1,\ldots, m_i-1\}$$ and let $R:\; \A^\N\to {\Bbb T}^k$ be the canonical coding map given by
$$
R(x)=\sum_{n=1}^\infty \Lambda^{-n}x_i,\quad x=(x_i)_{i=1}^\infty\in \A^\N.
$$
For any $\D\subseteq \A$, the set $R(\D^\N)$ is called a {\it self-affine Sierpinski sponge}. Whenever $k=2$, McMullen
\cite{McM84} and Bedford \cite{Bed84} determined the explicit value of the Hausdorff dimension of $R(\D^\N)$, and showed that there exists a Bernoulli product measure $\mu$  on $\D^\N$ such that $\dim_H \mu\circ R^{-1}=\dim_HR(\D^\N)$. Kenyon and Peres \cite{KePe96} extended this result to the general case $k\geq 2$, and
 moreover, they  proved  for each compact $T$-invariant set $K\subseteq {\Bbb T}^k$, there is an ergodic $\sigma$-invariant  $\mu$ on $\A^\N$ so that $\mu(R^{-1} (K))=1$ and
$\dim_H\mu\circ R^{-1}=\dim_HK$. Furthermore, Kenyon and Peres \cite{KePe96} proved the uniqueness of  $\mu\in \M(\D^\N,\sigma)$  satisfying $\dim_H \mu\circ R^{-1}=\dim_HR(\D^\N)$, by setting up  the following  formula for any ergodic $\eta\in \M(\A^\N,\sigma)$:
\begin{equation}
\label{e-KePe}
\dim_H\eta\circ R^{-1}=\frac{1}{\log m_k}h_\eta(\sigma)+\sum_{i=1}^{k-1}\left(\frac{1}{\log m_{k-i}}-\frac{1}{\log m_{k-i+1}}\right)h_{\eta\circ \tau_i^{-1}}(\sigma_{i}),
\end{equation}
where  $\tau_i$ denotes the one-block map from $\A^\N$ to $\A_i^\N$, with $\A_i=\prod_{j=1}^{k-i}\{0,1,\ldots, m_j-1\}$, so that each element in $\A$ (viewed as a $k$-dimensional vector)  is projected into to its first $(k-i)$ coordinates; and $\sigma_i$ denotes the left shift on $\A_i^\N$. Formula
\eqref{e-KePe} is an analogue of that for the Hausdorff dimension of $C^{1+\alpha}$ hyperbolic measures along unstable (respectively, stable) manifold established by Ledrappier and Young  \cite{LeYo85}.
 As Gatzouras and Peres pointed out in  \cite{GaPe96}, the uniqueness has not been known for more general invariant subsets $K$, even if $K=R(X)$, where $X\subset \A^\N$ is a general irreducible subshift of finite type. However, as a direct application of \eqref{e-KePe} and Theorem \ref{thm-1}, we have the following rather complete answer.
 \begin{thm}
\label{thm-2}
Let $X\subseteq \A^\N$ be a subshift satisfying  weak specification. Then there is a unique  $\mu\in \M(X,\sigma_X)$  such that $\dim_H\mu\circ R^{-1}=\dim_H R(X)$.
\end{thm}

Before this work, the problem of  Gatzouras and Peres had been studied and partially answered in the recent decade by different authors.    Assume that  $\pi$ is a factor map  between subshifts $X$ and $Y$, where $X$ is an irreducible subshift of finite type. Recall that  a {\it compensation function} for $\pi$ is a continuous function $F:\; X\to \R$ such that
$$
\sup_{\nu\in \M(Y,\sigma_Y)}\left(\nu(\phi)+h_\nu(\sigma_Y)\right)=\sup_{\mu\in \M(X,\sigma_X)}\left(\mu(\phi\circ \pi+F)+h_\mu(\sigma_X)\right)
$$
for all $\phi\in C(Y)$. Compensation functions were introduced in \cite{BoTu84} and studied systematically in \cite{Wal86}. Shin \cite{Shi01a} showed that if there exists a compensation function of the form $f\circ \pi$, with $f\in C(Y)$, and if $\frac{\alpha}{1+\alpha}f\circ \pi$ has a unique equilibrium state, then there is a unique measure $\mu$ maximizing the weighted entropy $h_\mu(\sigma_X)+\alpha h_{\mu\circ \pi^{-1}}(\sigma_Y)$. However, there exist
factor maps between irreducible subshifts  of finite type for which there are no such compensation functions \cite{Shi01b}.  Later, Petersen, Quas and Shin \cite{PQS03} proved that for each ergodic measure $\nu$ on $Y$, the number of ergodic measures $\mu$ of maximal entropy in the fibre $\pi^{-1}\{\nu\}$  is uniformly bounded; in particular,  if $\pi$ is a one-block map and  there is a symbol $b$ in the alphabet of $Y$ such that the pre-image of $b$ is a singleton (in this case,  $\pi: X\to Y$ is said to {\it have a singleton clump}), then
there is a unique  measures $\mu$ of maximal entropy in the fibre $\pi^{-1}\{\nu\}$ for each ergodic measure $\nu$ on $Y$. Recently, Yayama \cite{Yay09,Yay09a} showed the uniqueness of measures of  maximal weighted entropy if
 $\pi: X\to Y$  has a singleton clump. The uniqueness is further proved by Olivier \cite{Oli09} and Yayama
\cite{Yay09a} under an assumption that the projection of the ``Parry measure'' on $X$ has certain Gibbs property (however the assumption  only fulfils in some special cases).

The notions of weighted topological pressure and weighted equilibrium state were recently introduced by Barral and the author in \cite{BaFe09}, motivated from the study of the multifractal analysis on self-affine sponges \cite{Kin95,Ols98,BaMe07,BaMe08}. It was shown in \cite{BaFe09} that,  whenever $\pi_i:\; X_i\to X_{i+1}$ ($i=1,\ldots,k-1$) are one-block factor maps between one-sided full shifts $(X_i,\sigma_{X_i})$,  each $f\in V(\sigma_{X_1})$ has a unique $\ba$-weighted equilibrium state, which is Gibbs and mixing. The result had an interesting application  in the multifractal analysis \cite{BaFe09}. However, the approach given in \cite{BaFe09} depends upon  the simple fibre structure for the full shift case, and it does not work for the general case in Theorem \ref{thm-1}.

 The  main ingredient in our proof of Theorem \ref{thm-1} is to show the uniqueness of equilibrium states and conditional equilibrium states for certain sub-additive potentials, rather than for the classical additive potentials (or almost additive potentials).     A crucial step is to prove,  for certain functions $f$ defined on $\A^*$ (the set of finite words over $\A$), there exists an  ergodic  invariant  measures $\mu$ on the full shift space $\A^\N$ and $c>0$, so that $\mu(I)\geq c f(I)$ for $I\in \A^*$ (see Proposition \ref{pro-t1}).

The paper is organized as follows: In Sect.~\ref{S-2}, we introduce some basic notation and definitions about subshifts. In Sect.~\ref{S-3}, we present and prove some variational principles about certain sub-additive potentials. In Sect.~\ref{S-4}, we prove Proposition \ref{pro-t1}. In Sect.~\ref{S-5}, we prove the uniqueness of
equilibrium states for certain sub-additive potentials. In Sect.~\ref{S-6}, we prove the uniqueness of weighted equilibrium states for certain sub-additive potentials in the case $k=2$.  The extension to the general case $k\geq 2$ is given in Sect.~\ref{S-7}, together with the proof of Theorem \ref{thm-1}.

\section{Preliminaries about subshifts}
\label{S-2}

In this section, we  introduce some basic notation and definitions about subshifts.   The reader is referred to \cite{LiMa95} for the background and more details.

\subsection{One-sided subshifts over finite alphabets}
\label{S-2.1} Let $\A$ be a finite  set of symbols which we will call the {\it alphabet}.
 Let
$$
\A^*=\bigcup_{k=0}^\infty \A^k
$$
denote the set of all finite words with letters from $\A$, including the empty word $\varepsilon$. Let
$$
\A^\N=\left\{(x_i)_{i=1}^\infty:\; x_i\in \A\mbox { for }i\geq 1\right\}
$$
denote the collection of infinite sequences with entries from $\A$.
Then $\A^\N$ is a compact metric space endowed with the metric
$$
d(x,y)=2^{-\inf\{k:\; x_k\neq y_k\}},\quad x=(x_i)_{i=1}^\infty,\; y=(y_i)_{i=1}^\infty.
$$
For any $n\in \N$ and $I\in \A^n$, we write
\begin{equation}
\label{e-011}
[I]=\{(x_i)_{i=1}^\infty\in \A^\N:\; x_1\ldots x_n=I\}
\end{equation}
and call it an {\it $n$-th cylinder set} in $\A^\N$.

In this paper, a {\it topological dynamical system} is a continuous self map of a compact metrizable space.
The shift transformation $\sigma:\; \A^\N\to \A^\N$ is defined by $(\sigma x)_i=x_{i+1}$ for all $i\in \N$. The pair $(\A^\N,\sigma)$ forms a topological dynamical system which is called the {\it one-sided full shift over $\A$}.

If $X$ is a compact $\sigma$-invariant subset of $\A^\N$, that is, $\sigma(X)\subseteq X$, then the topological dynamical system $(X,\sigma)$ is called a {\it one-sided subshift over $\A$}, or simply, a {\it subshift}.  Sometimes, we denote a subshift $(X,\sigma)$ by $X$,  or $(X,\sigma_X)$.

A subshift $X$ over $\A$ is called a {\it subshift of finite type} if, there exists a matrix $A=(A(\alpha,\beta))_{\alpha,\beta\in \A}$ with entries $0$ or $1$ such that
$$X=\left\{(x_i)_{i=1}^\infty\in \A^\N:\; A(x_i,x_{i+1})=1\mbox{ for all } i\in \N\right\}.$$
If $A$ is irreducible (in the sense that, for any $\alpha, \beta\in \A$, there exists $n>0$ such that $A^n(\alpha,\beta)>0$), $X$ is called an {\it irreducible subshift of finite type}. Moreover if $A$ is primitive (in the sense that,  there exists $n>0$ such that $A^n(\alpha,\beta)>0$ for all $\alpha,\beta\in \A$), $X$ is called a {\it mixing subshift of finite type}.

The {\it language} $\lan (X)$ of a subshift $X$ is the set of all finite words (including the empty word $\varepsilon$) that occur as consecutive strings
$x_1\ldots x_n$ in the sequences $x=(x_i)_{i=1}^\infty$ which comprise $X$. That is,
$$
\lan(X)=\{I\in \A^*:\; I=x_1\ldots x_n\mbox{ for some  $x=(x_i)_{i=1}^\infty\in X$ and $n\geq 1$}\}\cup \{\varepsilon\}.
$$
Denote $|I|$ the length of a word $I$. For $n\geq 0$, denote
$$
\lan_n(X)=\{I\in \lan(X):\; |I|=n\}.
$$

Let $p\in \N$. A subshift $X$  is said to satisfy {\it $p$-specification} if for any $I, J\in \lan(X)$, there exists $K\in \lan_p(X)$  such that $IKJ\in \lan(X)$. We say that $X$ satisfies {\it specification} if it satisfies $p$-specification for some $p\in \N$.
  Similarly, $X$ is said to satisfy
{\it weak $p$-specification} if for any $I, J\in \lan(X)$, there exists $K\in \bigcup_{i=0}^p\lan_i(X)$  such that $IKJ\in \lan(X)$; and $X$ is said to satisfy {\it weak specification} if it satisfies weak $p$-specification for some $p\in \N$. It is easy to see that an irreducible subshift of finite type satisfies weak specification, whilst a  mixing subshift of finite type satisfies  specification.

Let $(X,\sigma_X)$ and  $(Y,\sigma_Y)$ be two subshifts over finite alphabets $\A$
and  $\A'$, respectively. We say that $Y$ is a {\it factor} of $X$
if,  there is a continuous surjective map $\pi:\; X\to Y$ such that
$\pi T=S\pi$.  Here $\pi$ is called a {\it factor map}. Furthermore
$\pi$ is called a {\it $1$-block  map} if there exists a map
$\pi:\; \A\to \A'$ such that
$$
\pi(x)=\left({\pi}(x_i)\right)_{i=1}^\infty,\qquad x=
(x_i)_{i=1}^\infty\in X.
$$
It is well known (see, e.g. \cite[Proposition 1.5.12]{LiMa95}) that  each factor map
$\pi:\; X\to Y$ between two subshifts $X$ and $Y$, will become a
$1$-block factor map if we enlarge the alphabet $\A$ and recode $X$
through a so-called {\it higher block representation of X}. Whenever $\pi:\; X\to Y$ is $1$-block, we write
$\pi I={\pi}(x_1)\ldots {\pi}(x_n)$ for
$I=x_1\ldots x_n\in \lan_n(X)$; clearly $\pi I\in \lan_n(Y)$.

\subsection{Two-sided subshifts over finite alphabets}
For a finite alphabet $\A$, let
$$
\A^\Z=\{x=(x_i)_{i\in \Z}:\; x_i\in \A \mbox{ for all }i\in \Z\}
$$
denote the collection of all bi-infinite sequence of symbols from $\A$.
Similarly, $\A^\Z$ is a compact metric space endowed with the metric
$$
d(x,y)=2^{-\inf\{|k|:\; x_k\neq y_k\}},\quad x=(x_i)_{i\in \Z},\; y=(y_i)_{i\in \Z}.
$$
The shift map $\sigma:\; \A^\Z\to \A^\Z$ is defined by $(\sigma x)_{i}=x_{i+1}$ for $x=(x_i)_{i\in \Z}$. The topological dynamical system $(\A^\Z,\sigma)$ is called the {\it two-sided full shift over $\A$}.

If $X\subseteq \A^\Z$ is compact and $\sigma(X)=X$, the topological dynamical system $(X,\sigma)$ is called a {\it two-sided subshift over $\A$}.

The definitions of $\lan(X)$, (weak) specification and factor maps for two-sided subshifts
 can be given in a way similar to the one-sided case.

\subsection{Some notation}
For two families of real numbers $\{a_i\}_{i\in \I}$ and
$\{b_i\}_{i\in \I}$,  we write
$$
\begin{array}{ll}
a_i\approx b_i \quad&\mbox{if there is $c>0$ such that $\frac{1}{c}
b_i\leq a_i\leq cb_i$ for $i\in \I$};\\
 a_i\succcurlyeq b_i \quad&\mbox{if there is $c>0$ such that $ a_i\geq
cb_i$ for $i\in
\I$};\\
 a_i\preccurlyeq b_i \quad&\mbox{if there is $c>0$ such that $ a_i\leq
cb_i$ for $i\in
\I$};\\
 a_i= b_i+O(1) \quad&\mbox{if there is $c>0$ such that
$|a_i-b_i|\leq c$ for $i\in \I$};\\
a_i\geq b_i+O(1) \quad&\mbox{if there is $c>0$ such that
$a_i-b_i\geq -c$ for $i\in \I$};\\
a_i\leq b_i+O(1) \quad&\mbox{if there is $c>0$ such that $a_i-b_i\leq c$
for $i\in \I$}.
\end{array}
$$

\section{Variational principles for sub-additive potentials}
\label{S-3}

In this section we present and prove some variational principles for certain sub-additive potentials. This is the starting point in our work.

First we give some notation and definitions.
Let $(X,\sigma_X)$ be a one-side subshift over a finite alphabet $\A$. We use $\M(X)$ to denote the set of
all Borel probability measures on $X$. Endow $\M(X)$  with the weak-star topology. Let   $\M(X,\sigma_X)$ denote  the set of
 all $\sigma_X$-invariant Borel probability measures on $X$. The sets $\M(X)$ and $\M(X,\sigma_X)$ are non-empty, convex and compact (cf. \cite{Wal82}).  For convenience, for $\mu\in \M(X)$ and $I\in \lan(X)$, we would like to write
 $$
 \mu(I):=\mu([I]\cap X),
 $$
where $[I]$ denotes the $n$-th cylinder in $\A^\N$ defined as in \eqref{e-011}.

 For $\mu\in \M(X,\sigma_X)$, the {\it measure theoretic entropy} of $\mu$ with respect to $\sigma_X$ is defined as
\begin{equation}
\label{e-v1}
 h_\mu(\sigma_X):=-\lim_{n\to \infty}\frac{1}{n}\sum_{I\in \lan_n(X)} \mu(I)\log \mu(I).
 \end{equation}
 The above limit exists since the sequence $(a_n)_{n=1}^\infty$, where $$a_n=-\sum_{I\in \lan_n(X)} \mu(I)\log \mu(I),$$
  satisfies $a_{n+m}\leq a_n+a_m$ for $n,m\in \N$. It follows
 \begin{equation}
 \label{e-v2}
 h_\mu(\sigma_X)=\inf_{n\in \N}\frac{1}{n}\sum_{I\in \lan_n(X)} \mu(I)\log \mu(I).
 \end{equation}
The function $\mu\mapsto h_\mu(\sigma_X)$ is affine and upper semi-continuous on $\M(X,\sigma_X)$ (cf. \cite{Wal82}).

 A sequence
$\Phi=(\log \phi_n)_{n=1}^\infty$ of functions on a subshift $X$ is called a {\it sub-additive
potential on $X$}, if  each $\phi_n$ is a non-negative continuous
function on $X$ and there exists $c>0$ such that
\begin{equation}
\label{e-em}
\phi_{n+m}(x)\leq c \phi_n(x)\phi_m(\sigma_X^nx),\quad \forall \; x\in X,\; n,m\in \N.
\end{equation}
For convenience, we denote by $\C_{sa}(X,\sigma_X)$ the collection of  sub-additive potentials on $X$. For $\Phi=(\log \phi_n)_{n=1}^\infty\in \C_{sa}(X,\sigma_X)$, define $\Phi_*:\; \M(X,\sigma_X)\to \R\cup \{-\infty\}$ by
\begin{equation}
\label{e-v3}
\Phi_*(\mu)=\lim_{n\to \infty}\frac{1}{n}\int \log \phi_n(x) d\mu(x).
\end{equation}
The limit in \eqref{e-v3} exists by the sub-additivity of $\log \phi_n$.
\begin{rem}{\rm
One observes that for $f\in C(X)$,
if $\Phi=(\log \phi_n)_{n=1}^\infty$ is given by $\phi_n(x)=\exp(S_nf(x))$, then $\Phi\in \C_{sa}(X,\sigma_X)$ and $\Phi_*(\mu)=\mu(f)$ for each $\mu\in \M(X,\sigma_X)$.}
\end{rem}

By the sub-additivity \eqref{e-em}, we have  the following  simple lemma (cf. Proposition 3.1 in \cite{FeHu09a}).
\begin{lem}
\label{lem-ba}
 \begin{itemize}
 \item[{\rm (i)}]
  $\Phi_*$ is affine and upper semi-continuous on $\M(X,\sigma_X)$.
  \item[{\rm (ii)}]
   There is a constant $C\in \R$ such that $\int \log \phi_n(x) d\mu(x) \geq n \Phi_*(\mu)-C$ for $n\in \N$ and $\mu\in \M(X,\sigma_X)$.
\end{itemize}
\end{lem}

\begin{de}
{\rm
 For $\Phi\in \C_{sa}(X,\sigma_X)$, $\mu\in \M(X,\sigma_X)$ is called an {\it equilibrium state of $\Phi$} if
$$
\Phi_*(\mu)+h_\mu(\sigma_X)=\sup\{\Phi_*(\eta)+h_\eta(\sigma_X):\; \eta\in \M(X,\sigma_X)\}.
$$
Let $\I(\Phi)$ denote the collection of all equilibrium states of $\Phi$.
}
\end{de}

A function $\phi:\; \lan(X)\to [0,\infty)$ is said to be {\it sub-multiplicative} if,  $\phi(\varepsilon)=1$ and there exists a constant $c>0$ such that $\phi(IJ)\leq c\phi(I)\phi(J)$ for any $IJ\in \lan(X)$. Furthermore, say
$\Phi=(\log \phi_n)_{n=1}^\infty\in \C_{sa}(X,\sigma_X)$ is {\it generated} by $\phi$ if
$$
\phi_n(x)=\phi(x_1\ldots x_n),\qquad x=(x_i)_{i=1}^\infty\in X.
$$

\begin{pro}
\label{pro-2.1}
Assume that  $\Phi=(\log \phi_n)_{n=1}^\infty\in \C_{sa}(X,\sigma_X)$ is generated by a sub-multiplicative
function $\phi:\; \lan(X)\to [0,\infty)$.  Then
\begin{itemize}
\item[{\rm (i)}] $\sup\{\Phi_*(\mu)+ h_\mu(\sigma_X):\; \mu\in \M(X,\sigma_X)\}=\lim_{n\to \infty} \frac1n\log u_n$, where $u_n$ is given by $$u_n=\sum_{I\in \lan_n(X)}\phi(I).$$

\item[{\rm (ii)}] $\I(\Phi)$ is a non-empty compact convex of $\M(X, \sigma_X)$. Furthermore each extreme point of  $\I(\Phi)$ is an ergodic measure.
\end{itemize}
\end{pro}

We remark that Proposition \ref{pro-2.1}(i) is a special case of Theorem 1.1 in \cite{CFH08} on the variational principle for sub-additive potentials. Proposition \ref{pro-2.1}(ii) actually holds for any $\Phi\in \C_{sa}(X,\sigma_X)$, by the affinity and upper semi-continuity of $\Phi_*(\cdot)$ and $h_{\cdot}(\sigma_X)$ on $\M(X,\sigma_X)$ (see the proof of Proposition \ref{pro-2.2}(ii) for details).

Now let $(X,\sigma_X)$ and $(Y,\sigma_Y)$ be  one-sided
subshifts over $\A, \A'$, respectively.  Assume that $Y$ is a factor of $X$ with a $1$-block factor map $\pi:\; X\to Y$.

\begin{de}
{\rm
For $\nu\in \M(Y,\sigma_Y)$, $\mu\in \M(X,\sigma_X)$ is called a {\it conditional equilibrium state of $\Phi$
 with respect to $\nu$}  if,  $\mu\circ \pi^{-1}=\nu$ and
$$
\Phi_*(\mu)+h_\mu(\sigma_X)=\sup\{\Phi_*(\eta)+h_\eta(\sigma_X):\; \eta\in \M(X,\sigma_X),\; \eta\circ \pi^{-1}=\nu\}.
$$
Let $\I_\nu(\Phi)$ denote the collection of all equilibrium states of $\Phi$ with respect to $\nu$.
}
\end{de}

The following result is a relativized version of Proposition \ref{pro-2.1}.

\begin{pro}
\label{pro-2.2}
Let $\Phi=(\log \phi_n)_{n=1}^\infty\in \C_{sa}(X,\sigma_X)$ be generated by a sub-multiplicative
function $\phi:\; \lan(X)\to [0,\infty)$. Let $\nu\in \M(Y, \sigma_Y)$. Then
\begin{itemize}
\item[{\rm (i)}] $\sup\{\Phi_*(\mu)+ h_\mu(\sigma_X)-h_\nu(\sigma_Y):\; \mu\in \M(X,\sigma_X),\;\mu\circ \pi^{-1}=\nu\}=\Psi_*(\nu)$, where $\Psi=(\log \psi_n)_{n=1}^\infty\in \C_{sa}(Y,\sigma_Y)$ is generated by
     a sub-multiplicative function $\psi:\; \lan(Y)\to [0,\infty)$, which satisfies
     \begin{equation}
     \label{e-em3}
     \psi(J)=\sum_{I\in \lan(X):\; \pi I=J}\phi(I),\; \quad J\in \lan(Y).
     \end{equation}

\item[{\rm (ii)}] $\I_\nu(\Phi)$ is a non-empty compact convex of $\M(X, \sigma_X)$. Furthermore, if $\nu$ is ergodic, then each extreme point of  $\I_\nu(\Phi)$ is an ergodic measure on $X$.
\end{itemize}
\end{pro}

We remark that Proposition \ref{pro-2.1} can be obtained from Proposition \ref{pro-2.2} by considering the special case that $Y$ is a singleton (correspondingly, $\A'$ consists of one symbol).

To prove Proposition \ref{pro-2.2}, we need the following lemmas.
\begin{lem}[\cite{Bow75}, p.~34]
\label{lem-b} Suppose $0\leq p_1,\ldots,p_m\leq 1$,
$s=p_1+\cdots+p_m\leq 1$ and $a_1,\ldots, a_m\geq 0$. Then
$$
\sum_{i=1}^m p_i(\log a_i-\log p_i)\leq s\log (a_1+\cdots
+a_m)-s\log s.
$$
\end{lem}

\begin{lem}[\cite{CFH08}, Lemma 2.3]
\label{lem-3}
Let $\Phi=(\log \phi_n)_{n=1}^\infty\in \C_{sa}(X,\sigma_X)$. Suppose $(\eta_n)_{n=1}^\infty$ is a sequence in
$\M(X)$. We form the new sequence $(\mu_n)_{n=1}^\infty$ by
$\mu_n=\frac 1n \sum_{i=0}^{n-1}\eta_n\circ \sigma_X^{-i}$. Assume that
$\mu_{n_i}$ converges to $\mu$ in $\M(X)$ for some subsequence
$(n_i)$ of natural numbers. Then $\mu\in \M(X,\sigma_X)$, and moreover
$$
\limsup_{i\to\infty}\frac{1}{n_i}\int \log
\phi_{n_i}(x)\;d\eta_{n_i}(x)\leq \Phi_*(\mu).
$$
\end{lem}

\begin{lem}[\cite{CFH08}, Lemma 2.4]
\label{lem-4} Denote $k=\#\A$. Then for any $\xi\in \M(X)$,  and  positive integers
$n,\ell$ with $n\geq 2\ell$, we have
$$
\frac 1n \sum_{I\in \lan_n(X)}\xi(I)\log \xi(I)\geq \frac
1\ell
\sum_{ I\in \lan_\ell(X)}\xi_n(I)\log \xi_n(I)-\frac{2\ell}{n}\log
k,
$$
where $\xi_n=\frac{1}{n}\sum_{i=0}^{n-1}\xi\circ \sigma_X^{-i}$.
\end{lem}

\bigskip
\noindent {\it Proof of Proposition \ref{pro-2.2}}. Fix $\nu\in \M(Y,\sigma_Y)$. For any $\mu\in \M(X,\sigma_X)$ with $\mu\circ \pi^{-1}=\nu$, and $n\in \N$,
we have
\begin{equation*}
\begin{split}
\sum_{I\in \lan_n(X)}&\mu(I)\log \phi(I)-\mu(I)\log \mu(I)\\
&=\sum_{J\in \lan_n(Y)}\sum_{I\in \lan_n(X):\; \pi I=J}\mu(I)\log \phi(I)-\mu(I)\log \mu(I)\\
&\leq \sum_{J\in \lan_n(Y)}\nu(J)\log \psi(J)-\nu(J)\log \nu(J)\quad\mbox{ (by Lemma \ref{lem-b})}.
  \end{split}
\end{equation*}
Dividing both sides by $n$ and letting $n\to \infty$,   we obtain
\begin{equation*}
\Phi_*(\mu)+h_\mu(\sigma_X)-h_\nu(\sigma_Y)\leq \Psi_*(\nu).
\end{equation*}
Thus to complete the proof of (i), it suffices to show that there exists  $\mu$ with $\mu\circ \pi^{-1}=\nu$, such that
$\Phi_*(\mu)+h_\mu(\sigma_X)-h_\nu(\sigma_Y)\geq \Psi_*(\nu)$. For this purpose, construct a sequence $(\eta_n)_{n=1}^\infty$  in $\M(X)$ such that
$$
\eta_n(I)=\frac{\nu(\pi I) \phi(I)}{\psi(\pi I)}, \quad \forall\; I\in \lan_n(X),
$$
where we take the convention $\frac{0}{0}=0$.  Clearly, $\eta_n\circ \pi^{-1}(J)=\nu(J)$ for all $J\in \lan_n(Y)$.
Set $\mu_n=\frac 1n \sum_{i=0}^{n-1}\eta_n\circ \sigma_X^{-i}$. Assume that
$\mu_{n_i}$ converges to $\mu$ in $\M(X)$ for some subsequence
$(n_i)$ of natural numbers. By Lemma \ref{lem-3}, $\mu\in \M(X,\sigma_X)$ and
\begin{equation}
\label{e-v4}
\begin{split}
\Phi_*(\mu)&\geq \limsup_{i\to\infty}\frac{1}{n_i}\int \log
\phi_{n_i}(x)\;d\eta_{n_i}(x)
=\limsup_{i\to\infty}\frac{1}{n_i}\sum_{I\in \lan_{n_i}(X)}\eta_{n_i}(I) \log
\phi(I).
\end{split}
\end{equation}

We first show that $\mu\circ \pi^{-1}=\nu$. Let $J\in \lan(Y)$. Denote $\ell=|J|$.
For  $n>\ell$ and $0\leq i\leq n-\ell$, we have
\[
\begin{split}
\eta_n\circ \sigma_X^{-i}\circ \pi^{-1}(J)&=\eta_n\circ \pi^{-1}\circ \sigma_Y^{-i}(J)\\
&=\sum_{J_1\in \lan_i(Y),\; J_2\in \lan_{n-i-\ell}(Y):\; J_1JJ_2\in \lan_n(Y)}\eta_n\circ \pi^{-1}(J_1JJ_2)\\
&=\sum_{J_1\in \lan_i(Y),\; J_2\in \lan_{n-i-\ell}(Y):\; J_1JJ_2\in \lan_n(Y)}\nu(J_1JJ_2)=\nu(J).
\end{split}
\]
It follows that $\mu_n\circ \pi^{-1}(J)=\frac 1n \sum_{i=0}^{n-1}\eta_n\circ \sigma_X^{-i}\circ \pi^{-1}(J)\to \nu(J)$, as $n\to \infty$. Therefore $\mu\circ \pi^{-1}(J)=\nu(J)$.  Since $J\in \lan(Y)$ is arbitrary, we have $\mu\circ \pi^{-1}=\nu$.

We next show that \begin{equation}
\label{e-v5}
\Phi_*(\mu)+h_\mu(\sigma_X)-h_\nu(\sigma_Y)\geq \Psi_*(\nu).
\end{equation}
Fix $\ell\in \N$. By Lemma \ref{lem-4}, we have for $n\geq 2\ell$,
$$
\frac 1n \sum_{I\in \lan_n(X)}\eta_n(I)\log \eta_n(I)\geq \frac
1\ell
\sum_{ I\in \lan_\ell(X)}\mu_n(I)\log \mu_n(I)-\frac{2\ell}{n}\log
k,
$$
where $k:=\#\A$. Since $\mu_{n_i}\to \mu$ as $i\to \infty$, we obtain
\begin{equation*}
\liminf_{i\to \infty}\frac{1}{n_i} \sum_{I\in \lan_{n_i}(X)}\eta_{n_i}(I)\log \eta_{n_i}(I)\geq \frac
1\ell
\sum_{ I\in \lan_\ell(X)}\mu(I)\log \mu(I).
\end{equation*}
Taking $\ell\to \infty$ yields
\begin{equation}
\label{e-v6}
\liminf_{i\to \infty}\frac{1}{n_i} \sum_{I\in \lan_{n_i}(X)}\eta_{n_i}(I)\log \eta_{n_i}(I)\geq -h_\mu(\sigma_X).
\end{equation}
Observe that
\begin{eqnarray*}
\sum_{I\in \lan_n(X)}\eta_n(I)\log \phi(I)&=&
\sum_{I\in \lan_n(X)}\eta_n(I)\log \frac{\eta_n(I)\psi(\pi I)}{\nu(\pi I)}\\
&=& \sum_{I\in \lan_n(X)}\eta_n(I)\log \eta_n(I)\\
&\mbox{}&\qquad \qquad  +\sum_{J\in \lan_n(Y)}\nu(J) (\log \psi(J)-\log \nu(J)).
\end{eqnarray*}
This together with \eqref{e-v6} yields
 $$
 \liminf_{i\to \infty}\frac{1}{n_i}\sum_{I\in \lan_{n_i}(X)}\eta_{n_i}(I)\log \phi(I)
 \geq -h_\mu(\sigma_X) +\Psi_*(\nu)+h_\nu(\sigma_Y).
$$
Applying \eqref{e-v4}, we have $\Phi_*(\mu)\geq -h_\mu(\sigma_X) +\Psi_*(\nu)+h_\nu(\sigma_Y)$. This proves \eqref{e-v5}. Hence the proof of (i) is complete.

Now we show (ii). By the above proof, we see that $\I_\nu(\Phi)\neq \emptyset$. The convexity of $\I_\nu(\Phi)$ follows directly from the affinity of $\Phi_*(\cdot)$ and $h_\cdot(\sigma_X)$ on $\M(X,\sigma_X)$. Furthermore, the compactness of $\I_\nu(\Phi)$ follows from the upper semi-continuity  of $\Phi_*(\cdot)$ and $h_\cdot(\sigma_X)$ on $\M(X,\sigma_X)$. Next,  assume that $\nu$ is ergodic and let $\mu$ be an extreme point of $\I_\nu(\Phi)$. We are going to show that $\mu$ is ergodic. Assume it is not true, that is, there exist $\mu_1, \mu_2\in \M(X,\sigma_X)$ with $\mu_1\neq \mu_2$, and $\alpha_1, \alpha_2\in (0,1)$ with $\alpha_1+\alpha_2=1$, such that $\mu=\sum_{i=1}^2\alpha_i \mu_i$.
Then $\nu=\mu\circ \pi^{-1}=\sum_{i=1}^2\alpha_i \mu_i\circ \pi^{-1}$. Since $\mu_i\circ \pi^{-1}\in \M(Y,\sigma_Y)$ for $i=1,2$ and $\nu$ is ergodic, we have  $\mu_1\circ \pi^{-1}=\mu_2\circ \pi^{-1}=\nu$.
Note that
$$
\Psi_*(\nu)=\Phi_*(\mu)+h_\mu(\sigma_X)-h_\nu(\sigma_Y)=
\sum_{i=1}^2 \alpha_i(\Phi_*(\mu_i)+h_{\mu_i}(\sigma_X)-h_\nu(\sigma_Y))
$$
and $\Phi_*(\mu_i)+h_{\mu_i}(\sigma_X)-h_\nu(\sigma_Y)\leq \Psi_*(\nu)$ by (i). Hence we have
$$\Phi_*(\mu_i)+h_{\mu_i}(\sigma_X)-h_\nu(\sigma_Y)=\Psi_*(\nu),\quad i=1,2.
$$
That is, $\mu_i\in \I_\nu(\Phi)$ for $i=1,2$. However $\mu=\sum_{i=1}^2\alpha_i\mu_i$. It contradicts the assumption that $\mu$ is an extreme point of $\I_\nu(\Phi)$.
This finishes the proof of the proposition.
\qed

\begin{de}
{\rm
Let $\ba=(a_1,a_2)\in \R^2$ with $a_1>0$ and $a_2\geq 0$. For $\Phi\in \C_{sa}(X,\sigma_X)$,
 $\mu\in \M(X,\sigma_X)$ is called an {\it ${\ba}$-weighted equilibrium state of $\Phi$ for the factor map $\pi$}, or simply, {\it ${\ba}$-weighted equilibrium state of $\Phi$}, if
\begin{equation}
\label{e-1.9}
\begin{split}
\Phi_*(\mu)&+a_1h_\mu(\sigma_X)+a_2h_{\mu\circ \pi^{-1}}(\sigma_Y)\\
&=\sup\{\Phi_*(\eta)+a_1h_\eta(\sigma_X)+a_2h_{\eta\circ \pi^{-1}}(\sigma_Y):\; \eta\in \M(X,\sigma_X)\}.
\end{split}
\end{equation}
We use
$\I(\Phi,\ba)$ to denote the collection of all $\ba$-weighted equilibrium states of $\Phi$.
The value in the right hand side of \eqref{e-1.9} is called the {\it $\ba$-weighted topological pressure of $\Phi$} and is denoted by $P^\ba(\sigma_X,\Phi)$.
Each $\mu\in \I(\Phi,\ba)$ is called an {\it {\ba}-weighted equilibrium state} of $\Phi$.
}
\end{de}
As a corollary of Propositions \ref{pro-2.1} and \ref{pro-2.2}, we have
\begin{cor}
\label{cor-1.1}
Let $\Phi=(\log \phi_n)_{n=1}^\infty\in \C_{sa}(X,\sigma_X)$ be generated by a sub-multiplicative
function $\phi:\; \lan(X)\to [0,\infty)$. Define $\phi^{(2)}:\; \lan(Y)\to [0,\infty)$  by
\begin{equation*}
\label{e-}
\phi^{(2)}(J)=\Big(\sum_{I\in \lan_n(X):\; \pi I=J}\phi(I)^{\frac{1}{a_1}}\Big)^{a_1} \mbox{ for } J\in \lan_n(Y),\; n\in \N.
\end{equation*}
Let $\Phi^{(2)}=(\log \psi_n)_{n=1}^\infty\in \C_{sa}(Y,\sigma_Y)$ be generated by $\phi^{(2)}$. Then
\begin{itemize}
\item[{\rm (i)}]
$\mu\in \I(\Phi,\ba)$ if and only if  $\mu\circ \pi^{-1}\in \I(\frac{1}{a_1+a_2}\Phi^{(2)})$ and
$\mu\in \I_{\mu\circ \pi^{-1}}(\frac{1}{a_1}\Phi)$, where
$\frac{1}{a_1+a_2}\Phi^{(2)}:=(\log (\psi_n^{1/(a_1+a_2)}))_{n=1}^\infty$ and
$\frac{1}{a_1}\Phi:=(\log (\phi_n^{1/a_1}))_{n=1}^\infty$.
\item[{\rm (ii)}]
Furthermore, $\I(\Phi,\ba)$ is a non-empty compact convex set, and each extreme point of $\I(\Phi,\ba)$ is ergodic.
\item[{\rm (iii)}]  $\I(\Phi,\ba)$ is a singleton if and only if $\I(\frac{1}{a_1+a_2}\Phi^{(2)})$ is a singleton $\{\nu\}$ and, $\I_{\nu}(\frac{1}{a_1}\Phi)$ contains a unique ergodic measure.
\end{itemize}
\end{cor}
\Proof
Note that for $\mu\in \M(X,\sigma_X)$,
\[
\begin{split}
\Phi_*(\mu)&+a_1h_\mu(\sigma_X)+a_2h_{\mu\circ \pi^{-1}}(\sigma_Y)\\
&=
\Phi_*(\mu)+a_1(h_\mu(\sigma_X)-h_{\mu\circ \pi^{-1}}(\sigma_Y))+
(a_1+a_2)h_{\mu\circ \pi^{-1}}(\sigma_Y).
\end{split}
\]
By Proposition \ref{pro-2.2},
\begin{equation*}
\begin{split}
\sup &
\left
\{
\Phi_*(\eta)+a_1(h_\eta(\sigma_X)-h_{\eta\circ \pi^{-1}}(\sigma_Y)):\;\eta\in \M(X,\sigma_X),\; \eta\circ \pi^{-1}=\mu\circ \pi^{-1}
\right
\}
\\
&=\Phi^{(2)}_*(\mu\circ \pi^{-1}).
\end{split}
\end{equation*}
Hence $\mu\in \I(\Phi,\ba)$ if and only if that (i)
$$\Phi_*(\mu)+a_1(h_\mu(\sigma_X)-h_{\mu\circ \pi^{-1}}(\sigma_Y))
=\Phi^{(2)}_*(\mu\circ \pi^{-1});$$ and (ii) $$\Phi^{(2)}_*(\mu\circ \pi^{-1})+(a_1+a_2)h_{\mu\circ \pi^{-1}}(\sigma_Y)
=\sup_{\nu\in \M(Y,\sigma_Y)}\left(\Phi^{(2)}(\nu)+(a_1+a_2)h_\nu(\sigma_Y)\right).$$
That is, $\mu\in \I(\Phi,\ba)$ if and only if $\mu\in \I_{\mu\circ \pi^{-1}}(\frac{1}{a_1}\Phi)$ and $\mu\circ \pi^{-1}\in \I(\frac{1}{a_1+a_2}\Phi^{(2)})$. This proves (i). The proof of (ii) is essentially identical  to that of Proposition \ref{pro-2.1}(ii). Part (iii) follows from (i) and (ii).
\qed

\begin{rem}
{\rm
Proposition \ref{pro-2.2} was proved in \cite{BaFe09} in the special case that $\pi:\; X\to Y$ is a one-block factor map between full shifts.
Independently, Proposition \ref{pro-2.2} and Corollary \ref{cor-1.1}  were set up in \cite{Yay09a} for the special case that $\phi\equiv 1$ and $X$ is an irreducible subshift of finite type, by a direct combination of \cite[Theorem 2.1]{LeWa77} and  \cite[Corollary]{PeSh05}.

}
\end{rem}

\section{Ergodic invariant measures associated with certain functions on $\A^*$}
\label{S-4}
Let $\A$ be a finite alphabet and let  $\A^*=\bigcup_{n=0}^\infty \A^n$.
We define two collections of functions over $\A^*$.

\begin{de}
\label{de-1}
Let $p\in \N$. Define  $\Omega_w(\A^*, p)$ to be the collection of
functions  $f: \A^*\to [0,1]$ such that there exists $c>0$ so that
\begin{itemize}
\item[{\rm (H1)}] $\sum_{I\in \A^n}f(I)=1$ for any $n\geq 0$.
\item[{\rm (H2)}] For any $I,J\in \A^*$, there exists $K\in \bigcup_{i=0}^p\A^i$ such that  $f(IKJ)\geq cf(I)f(J)$.
\item[{\rm (H3)}] For each $I\in \A^*$, there exist $i,j\in \A$ such that
$$f(iI)\geq cf(I),\quad f(Ij)\geq cf(I).$$
\end{itemize}
\end{de}

\begin{de}
Let $p\in \N$. Let $\Omega(\A^*, p)$ denote the collection of
functions $g: \A^*\to [0,1]$ such that there exists $c>0$ so that
\begin{itemize}
\item[{\rm (A1)}]   $\sum_{I\in \A^n}g(I)=1$ for any $n\geq 0$.
\item[{\rm (A2)}] For any $I,J\in \A^*$, there exists $K\in \A^p$ such that  $g(IKJ)\geq cg(I)g(J)$.
\end{itemize}
\end{de}

 For $f\in \Omega_w(\A^*, p)\cup \Omega(\A^*, p)$, define a map $f^*:\; \A^*\to
 [0,\infty)$ by
\begin{equation}
\label{e-t1}
f^*(I)=\sup_{m,n\geq 0}f_{m,n}(I), \qquad I\in \A^*,
\end{equation}
where $f_{m,n}(I):=\sum_{I_1\in \A^m}\sum_{I_2\in \A^n} f(I_1II_2)$.
Clearly,  $f(I)=f_{0,0}(I)\leq f^*(I)\leq 1$ for any $I\in \A^*$.

The main result in this section is the following  proposition, which plays a key role in our proof of Theorem \ref{thm-2}.
\begin{pro}
\label{pro-t1} Let $f\in \Omega_w(\A^*, p)\cup \Omega(\A^*, p)$ and $f^*$ be defined as in \eqref{e-t1}. Let
$(\eta_n)_{n=1}^\infty$ be a sequence of Borel probability measures
on $\A^\N$ satisfying
$$
\eta_n(I)=f(I),\qquad\forall\;  I\in \A^n.
$$
We form the new sequence $(\mu_n)_{n=1}^\infty$ by
$\mu_n=\frac{1}{n}\sum_{i=0}^{n-1}\eta_n\circ \sigma^{-n}$. Assume that $\mu_{n_i}$ converges to $\mu$ for some subsequence $(n_i)$ of natural numbers.  Then  $\mu\in \M(\A^\N,\sigma)$ and it satisfies the following
properties:
\begin{itemize}
\item[{\rm (i)}]There is a  constant $C_1> 0$ such that $C_1f^*(I)\leq \mu(I)\leq  f^*(I)$ for all $I\in \A^*$.
\item[{\rm (ii)}] There is a constant $C_2>0$ such that
$$\liminf_{n\to \infty}\sum_{i=0}^p\mu(A\cap \sigma^{-n-i}(B))\geq C_2\mu(A)\mu(B)$$ for any Borel sets $A,B\subseteq \A^\N$.
\item[{\rm (iii)}] $\mu$ is ergodic.
\item[{\rm (iv)}] $\mu$ is the unique ergodic  measure on $\A^\N$ such that  $\mu(I)\geq C_3 f(I)$ for all $I\in \A^*$
and  some constant $C_3>0$.
\item[{\rm (v)}] $\frac{1}{n}\sum_{i=0}^{n-1}\eta_n\circ \sigma^{-n}$
converges to $\mu$, as $n\to \infty$.
\end{itemize}
Furthermore if  $f\in \Omega(\A^*, p)$, we have
 \begin{itemize}
\item[{\rm (vi)}]  There is a constant $C_4>0$ such that
$$\liminf_{n\to \infty}\mu(A\cap \sigma^{-n}(B))\geq C_4\mu(A)\mu(B)$$ for any Borel sets $A,B\subseteq \A^\N$.
\end{itemize}
\end{pro}

To prove the above proposition, we need several lemmas.
\begin{lem}
\label{lem-t1}
 Let $f\in \Omega_w(\A^*, p)\cup \Omega(\A^*, p)$. Then there is a constant $C>0$, which depends on $f$, such that
 \begin{itemize}
 \item[{\rm (i)}] $f_{m',n'}(I)\geq C f_{m,n}(I)$ for any $I\in \A^*$, $m'\geq m+p$ and $n'\geq n+p$.
 \item[{\rm (ii)}] For each $I\in \A^*$, there exists an integer $N=N(I)$ such that
$$
f_{m,n}(I)\geq (C/2) f^*(I),\qquad \forall \; m,n\geq N.
$$
\end{itemize}
\end{lem}

\Proof
   To show (i), we first assume $f\in \Omega_w(\A^*, p)$. Let $c$ be the constant associated with $f$ in Definition
   \ref{de-1}.
    Fix $I\in \A^*$ and $m,n,m',n'\in \N\cup\{0\}$ such that $m'\geq m+p$ and $n'\geq n+p$. By (H2), for given $I_1\in \A^m$, $I_2\in \A^n$, $I_3\in \A^{m'-m-p}$ and $I_4\in \A^{n'-n-p}$,  there exist $K_1,K_2\in \bigcup_{i=0}^p \A^i$ so that $$
f(I_3K_1I_1II_2 K_2I_4)\geq c^{2}f(I_3)f(I_1II_2)f(I_4).
$$
Furthermore by (H3), there exist $K_3,K_4\in \bigcup_{i=0}^p\A^i$ so that $|K_1|+|K_3|=p$,  $|K_2|+|K_4|=p$ and
\begin{equation}
\label{e-p1}
f(K_3I_3K_1I_1II_2 K_2I_4K_4)\geq c^{2p}f(I_3K_1I_1II_2 K_2I_4)\geq c^{2p+2} f(I_3)f(I_1II_2)f(I_4).
\end{equation}
Summing over $I_1\in \A^m$, $I_2\in \A^n$, $I_3\in \A^{m'-m-p}$ and $I_4\in \A^{n'-n-p}$,  and using (H1), we obtain
$$
f_{m',n'}(I)\geq \frac{1}{M}c^{2p+2}f_{m,n}(I),
$$
where $M$ denotes the number of different tuples $(J_1,J_2,
J_3,J_4)\in (\A^*)^4$ with $|J_1|+|J_3|=p$ and $|J_2|+|J_4|=p$.

Now assume $f\in \Omega(\A^*, p)$. Instead of (\ref{e-p1}),  by
(A2), we can find $K_1,K_2\in \A^p$ such that
\begin{equation*}
f(I_3K_1I_1II_2 K_2I_4)\geq c^{2}f(I_3)f(I_1II_2)f(I_4).
\end{equation*}
Summing over $I_1, I_2, I_3, I_4$  yields
$$
f_{m',n'}(I)\geq c^{2}f_{m,n}(I).
$$
This proves (i) by taking $C=\min\{c^2,
\frac{1}{M}c^{2p+2}\}=\frac{1}{M}c^{2p+2}$.

To show (ii), note that  $f^*(I)=\sup_{m,n\geq 0}f_{m,n}(I)$. Hence we can pick $m_0,n_0$ such that $f_{m_0,n_0}(I)\geq f^*(I)/2$. Let $N=m_0+n_0+p$. Then by (i), for any $m,n\geq N$, we have
$$f_{m,n}(I)\geq C f_{m_0,n_0}(I)\geq \frac{C}{2}f^*(I).$$
This finishes the proof of prove the lemma.
\qed

\begin{lem}
\label{lem-t2} Let $f\in \Omega_w(\A^*, p)\cup \Omega(\A^*, p)$.
Then there exists a constant $C'>0$ such that for any $I,J\in \A^*$,
there exists an integer $N=N(I,J)$ such that
$$
\sum_{i=0}^p \sum_{K\in \A^{n+i}}f^*(IKJ)\geq C'f^*(I)f^*(J),\qquad \forall \; n\geq N.
$$
In particular,  if $f\in \Omega(\A^*, p)$, then the above inequality
can be strengthened as
$$
\sum_{K\in \A^{n}}f^*(IKJ)\geq C'f^*(I)f^*(J),\qquad \forall \; n\geq N.
$$

\end{lem}
\Proof
First assume $f\in \Omega_w(\A^*, p)$. Let $C$ be the constant associated with $f$ in Lemma
\ref{lem-t1}.  Fix $I,J\in \A^*$. By Lemma \ref{lem-t1}(ii),  there
exists $k\in \N$ such that for  $m_1,m_2,m_3, m_4\geq k$,
$$f_{m_1,m_2}(I)\geq \frac{C}{2}f^*(I),\quad  f_{m_3,m_4}(J)\geq \frac{C}{2}f^*(J).$$
Take $N=2k$. Let $n\geq N$. Then we have
$$
f_{k,n-k}(I)\geq \frac{C}{2}f^*(I), \quad f_{k,k}(J)\geq \frac{C}{2}f^*(J).
$$
By (H2),  for any $I_1\in \A^k$, $I_2\in \A^{n-k}$, $J_1,J_2\in \A^k$,
we have
\begin{equation}
\label{e-p2}
\sum_{i=0}^p\sum_{U\in \A^i}f(I_1II_2UJ_1JJ_2)\geq cf(I_1II_2)f(J_1JJ_2).
\end{equation}
Summing over $I_1,I_2,J_1,J_2$ yields
$$\sum_{i=0}^p\sum_{K\in \A^{n+i}}f_{k,k}(IKJ)\geq cf_{k,n-k}(I)f_{k,k}(J).$$
Hence, we have
$$
\sum_{i=0}^p\sum_{K\in \A^{n+i}}f^*(IKJ)\geq
cf_{k,n-k}(I)f_{k,k}(J)\geq c(C/2)^2f^*(I)f^*(J).
$$

Next assume $f\in \Omega(\A^*, p)$.  By (A2), instead of
(\ref{e-p2}), we have
 $$
 \sum_{U\in \A^p}f(I_1II_2UJ_1JJ_2)\geq cf(I_1II_2)f(J_1JJ_2)
 $$
 for any $I_1\in \A^k$, $I_2\in \A^{n-k}$, $J_1,J_2\in \A^k$.
 Summing over $I_1,I_2,J_1,J_2$ we obtain
$$\sum_{K\in \A^{n+p}}f_{k,k}(IKJ)\geq cf_{k,n-k}(I)f_{k,k}(J)\geq c(C/2)^2f^*(I)f^*(J).$$
Hence $\sum_{K\in \A^{n+p}}f^*(IKJ)\geq  c(C/2)^2f^*(I)f^*(J)$. This finishes the proof of the lemma.
\qed

\bigskip
\noindent{\it Proof of Proposition \ref{pro-t1}}.
 By \cite[Theorem 6.9]{Wal82}, $\mu$ is $\sigma$-invariant.
 Fix $I\in \A^*$.  Let $m=|I|$. For $n>m$,  we have
 \begin{equation*}
 \begin{split}
 \mu_n(I)&=\frac{1}{n}\left(\sum_{i=0}^{n-m}\eta_n\circ \sigma^{-i}(I)+\sum_{j=n-m+1}^{n-1}\eta_n\circ\sigma^{-j}(I)\right)\\
&=\frac{1}{n}\left(\sum_{i=0}^{n-m}f_{i,n-m-i}(I)+\sum_{j=n-m+1}^{n-1}\eta_n\circ\sigma^{-j}(I)\right).
 \end{split}
 \end{equation*}
 Applying Lemma \ref{lem-t1}(ii) to the above equality yields
 $$
 \frac{C}{2} f^*(I) \leq \liminf_{n\to \infty}\mu_n(I)\leq \limsup_{n\to \infty}\mu_n(I)\leq f^*(I),
 $$
 where $C>0$ is a constant independent of $I$. Hence
 $$(C/2)f^*(I)\leq \mu(I)\leq f^*(I).$$
  This proves (i) by taking $C_1=C/2$.

   By (i) and Lemma  \ref{lem-t2}, we have
 \begin{equation}
 \label{e-m1}
 \begin{split}
 \liminf_{n\to \infty} &
 \sum_{i=0}^p\
 \mu([I]\cap \sigma^{-n-i}([J])) \geq
 C_1
 \liminf_{n\to \infty}\sum_{i=0}^p\sum_{K\in \A^{n+i}} f^*(IKJ)\\
 &\geq  C_1C'f^*(I)f^*(J)\geq  C_1C' \mu(I)\mu(J)
 \end{split}
 \end{equation}
 for some constant $C'>0$ and all $I,J\in \A^*$. Take $C_2=C_1C'$.
 Since $\{[I]:\; I\in \A^*\}$ generates the Borel $\sigma$-algebra of $\A^\N$,   (ii) follows from (\ref{e-m1})
 by a standard argument.

 As a consequence of (ii), for any Borel sets $A,B\subseteq \A^\N$ with $\mu(A)>0$ and $\mu(B)>0$, there exists $n$ such that
 $\mu(A\cap \sigma^{-n}(B))>0$. This implies that $\mu$ is ergodic (cf. \cite[Theorem 1.5]{Wal82}). This proves (iii).

 To prove (iv), assume that $\eta$ is an ergodic  measure on $\A^\N$ so
 that there exists $C_3>0$ such that
 $$\eta(I)\geq C_3 f(I),\quad \forall\; I\in \A^*.$$
 Then for any $I\in \A^*$ and $m,n\in \N$,
 $$
 \eta(I)= \sum_{I_1\in \A^m}\sum_{I_2\in \A^n}\eta(I_1II_2)
 \geq C_3\sum_{I_1\in \A^m}\sum_{I_2\in
 \A^n}f(I_1II_2)=C_3f_{m,n}(I).
 $$
 Hence $\eta(I)\geq C_3f^*(I)\geq C_3\mu(I)$. It implies that $\mu$ is absolutely continuous with respect to $\eta$.
 Since any two different ergodic measures on $\A^\N$ are
singular to each other (cf. \cite[Theorem 6.10(iv)]{Wal82}), we have $\eta=\mu$. This proves
(iv). Notice that (v) follows directly from (i), (iii) and (iv).

  Now assume that $f\in \Omega(\A^*, p)$. Instead of (\ref{e-m1}), by (i) and Lemma \ref{lem-t2} we have
\begin{equation*}
 \begin{split}
\liminf_{n\to \infty}~&
\mu([I]\cap\sigma^{-n}([J]))=C_1\liminf_{n\to \infty}\sum_{K\in \A^{n}} f^*(IKJ)\\
 &\geq  C_1C'f^*(I)f^*(J)\geq  C_1C'
 \mu(I)\mu(J)=C_2\mu(I)\mu(J),
 \end{split}
\end{equation*}
from which (vi) follows.   This finishes the proof of Proposition \ref{pro-t1}.
  \qed

\section{Equilibrium states for certain sub-additive potentials}
\label{S-5}

In this section, we show the uniqueness of equilibrium states for certain sub-additive potentials on one-sided subshifts.

Let $(X,\sigma_X)$ be a subshift over a finite alphabet $\A$.  For $n\geq 1$, denote
$$
\lan_n(X)=\{I\in \A^n:\; X\cap [I]\neq \emptyset\}.
$$
Denote $\lan_0(X)=\{\varepsilon\}$, where $\varepsilon$ denotes the empty word.  Set $\lan(X)=\bigcup_{i=0}^\infty \lan_n(X)$.

Let $p\in \N$. We use $\D_w(X, p)$ denote the collection of functions  $\phi:\; \lan(X)\to [0,\infty)$ such that
 $\phi(I)>0$ for at least one $I\in \lan(X)\backslash \{\varepsilon\}$, and there exist $0<c\leq 1$ so that
\begin{itemize}
\item[(1)] $\phi(IJ)\leq c^{-1}\phi(I)\phi(J)$
for any $IJ\in \lan(X)$.
\item[(2)]  For any $I,J\in \lan(X)$, there exists $K\in \bigcup_{i=0}^p\lan_i(X)$ such that $IKJ\in \lan(X)$ and $\phi(IKJ)\geq c\phi(I)\phi(J)$.
\end{itemize}

Furthermore, we use $\D(X, p)$ denote the collection of functions  $\phi:\; \lan(X)\to [0,\infty)$ such that
 $\phi(I)>0$ for at least one $I\in \lan(X)\backslash \{\varepsilon\}$, and there exist $0<c\leq 1$ so that
$\phi$ satisfies the above condition (1), and
\begin{itemize}
\item[(2')]  For any $I,J\in \lan(X)$, there exists $K\in \lan_p(X)$ such that $IKJ\in \lan(X)$ and $\phi(IKJ)\geq c\phi(I)\phi(J)$.
\end{itemize}
\begin{rem}
{\rm \begin{itemize}\item[(i)] $\D(X,p)\subseteq \D_w(X,p)$.
 \item[(ii)]$\D_w(X,p)\neq \emptyset$ if and only if $X$ satisfies weak $p$-specification. The necessity is obvious. For   the sufficiency,  if $X$ satisfies weak $p$-specification, then the constant function $\phi\equiv 1$ on $\lan(X)$
is an element in $\D_w(X,p)$. Similarly, $\D(X,p)\neq \emptyset$ if and only if $X$ satisfies $p$-specification.
\end{itemize}
}
\end{rem}

\begin{lem}
\label{lem-g1}
Suppose $\phi\in \D_w(X,p)$. Then the following two properties hold:
\begin{itemize}
\item[{\rm (i)}]  There exists a constant $\gamma>0$ such that for each $I\in \lan(X)$, there exist $i,j\in \A$ such that $\phi(iI)\geq \gamma \phi(I)$ and $\phi(Ij)\geq \gamma\phi(I)$.
\item[{\rm (ii)}]
Let $u_n=\sum_{J\in X_n}\phi(J)$. Then  the limit $u=\lim_{n\to \infty}
(1/n)\log u_n$ exists and  $u_{n}\approx \exp(nu)$.
\end{itemize}
\end{lem}
\Proof Let $\phi\in \D_w(X,p)$ with the corresponding constant $c\in (0,1]$.  For (i), we only prove there exists a constant $\gamma>0$ such that for each $I\in \lan(X)$, there exist $j\in \A$ such that $\phi(Ij)\geq \gamma \phi(I)$. The other statement (there exists  $i\in \A$ so that $\phi(iI)\geq \gamma \phi(I)$)  follows by an identical argument.
Fix a word $W\in \lan(X)\backslash \{\varepsilon\}$ such that $\phi(W)>0$. Let $I\in \lan(X)$ so that $\phi(I)>0$. Then there exists $K\in \bigcup_{i=0}^p\lan_i(X)$ such that $\phi(IKW)\geq c\phi(I)\phi(W)$. Write $KW=jU$, where $j$ is the first letter in the word $KW$. Then $$\phi(Ij)\phi(U)\geq c\phi(IjU)=c\phi(IKW)\geq c^2\phi(I)\phi(W).$$ Hence $\phi(U)>0$ and $\phi(Ij)\geq c^2\phi(I)\phi(W)/\phi(U)$. Since there are only finite  possible $U$ (for $|U|\leq |W|+p$), $\phi(Ij)/\phi(I)\geq \gamma$ for some constat $\gamma>0$.

To see (ii), we have
\begin{equation}
\label{e-aa1}
\begin{split}
u_{n+m}&=\sum_{I\in \lan_n(X), \;J\in \lan_m(X):\; IJ\in \lan_{n+m}(X)}
\phi(IJ)\\
&\leq \sum_{I\in \lan_n(X), \;J\in \lan_m(X)} c^{-1}\phi(I)\phi(J)
=c^{-1}u_nu_m
\end{split}
\end{equation}
and
\begin{equation}
\label{e-aa2}
\begin{split}
\sum_{k=0}^pu_{n+m+k}&=\sum_{I\in \lan_n(X), \;J\in \lan_m(X)}\; \sum_{K\in \bigcup_{i=0}^p\lan_i(X):\;  IKJ\in \lan(X)}\phi(IKJ)\\
&\geq
\sum_{I\in \lan_n(X),\; J\in \lan_m(X)} c\phi(I)\phi(J)=cu_nu_m.
\end{split}
\end{equation}
 On the other hand,   $$u_{n+1}=\sum_{I\in \lan_n(X)}\;\sum_{j\in \A:\; Ij\in \lan_{n+1}(X)}\phi(Ij)\geq \gamma \sum_{I\in \lan_n(X)}\phi(I)=\gamma u_n,$$
  and $u_{n+1}\leq c^{-1}u_1u_n$ by \eqref{e-aa1}. Hence $u_{n+1}\approx u_{n}$. This together with
 \eqref{e-aa1} and \eqref{e-aa2} yields $u_{n+m}\approx u_nu_m$, from which (ii) follows.
\qed

Note that we have introduced $\Omega_w(\A^*,p)$ and $\Omega(\A^*,p)$ in Sect.~\ref{S-4}.
As a direct consequence of Lemma \ref{lem-g1}, we have
\begin{lem}
\label{lem-dir}
Let $\phi\in \D_w(X,p)$. Define
$f:\A^*\to [0,1]$ by
\begin{equation}
\label{e-aa3}
f(I)=\left\{
\begin{array}{cl}
\displaystyle\frac{\phi(I)}{\sum_{J\in \lan_n(X)}\phi(J)} &\mbox{ if }I\in \lan_n(X),\; n\geq 0,\\
&\\
0&\mbox{ if }I\in \A^*\backslash \lan(X).
\end{array}
\right.
\end{equation}
Then $f\in \Omega_w(\A^*,p)$, and $f(IJ)\preccurlyeq f(I)f(J)$ for $I, J\in \A^*$.
Moreover if $\phi\in \D(X,p)$, then $f\in \Omega(\A^*,p)$.
\end{lem}

\begin{lem}
\label{lem-4.1}
Let $\eta,\mu\in \M(X,\sigma_X)$. Assume that $\eta$ is not absolutely continuous with respect to $\mu$. Then
$$
\lim_{n\to \infty}\sum_{I\in \lan_n(X)} \eta(I) \log \mu(I)-\eta(I)\log \eta(I)=-\infty.
$$
\end{lem}
\Proof
We take a slight modification of the proof of Theorem 1.22 in \cite{Bow75}.
Since $\eta$ is not absolutely continuous with respect to $\mu$, there exists $c\in (0,1)$ such that for any $0<\epsilon<c/2$, there exists a Borel set $A\subset X$ so that
$$
\eta(A)>c\quad  \mbox{and} \quad \mu(A)<\epsilon.
$$
Applying \cite[Lemma 1.23]{Bow75}, we see that for each sufficiently large $n$, there exists $F_n\subset \lan_n(X)$ so that
$$
\mu(A\vartriangle A_n)+\eta(A\vartriangle A_n)<\epsilon\quad
\mbox{with }A_n:=\bigcup_{I\in F_n}[I]\cap X,
$$
which implies $\eta(A_n)>c-\epsilon>c/2$ and $\mu(A_n)<2\epsilon$. Using Lemma \ref{lem-b}, we obtain
\begin{equation}
\label{e-t3.1}
\begin{split}
\sum_{I\in F_n}\eta(I) \log \mu(I)-\eta(I)\log \eta(I)
&\leq \eta(A_n)\log \mu(A_n)+\sup_{0\leq s\leq 1}s\log (1/s)\\
& \leq (c/2)\log (2\epsilon)+\log 2
\end{split}
\end{equation}
and
\begin{equation}
\label{e-t3.2}
\begin{split}
&\sum_{I\in \lan_n(X)\backslash F_n}\eta(I) \log \mu(I)-\eta(I)\log \eta(I)\\
&\quad \leq \eta(X\backslash A_n)\log \mu(X\backslash A_n)+\sup_{0\leq s\leq 1}s\log (1/s)\leq \log 2.
\end{split}
\end{equation}
Combining (\ref{e-t3.1}) and (\ref{e-t3.2})  yields
$$
\sum_{I\in \lan_n(X)}\eta(I) \log \mu(I)-\eta(I)\log \eta(I) \leq (c/2)\log (2\epsilon)+2\log 2,
$$
from which the lemma follows.
\qed

The main result in this section is the following

\begin{thm}
 \label{thm-3.1}
 Let $\phi\in \D_w(X, p)$.  Let $\Phi=(\log \phi_n)_{n=1}^\infty\in \C_{sa}(X,\sigma_X)$ be generated by $\phi$, i.e. $\phi_n(x)=\phi(x_1\ldots x_n)$ for $x=(x_i)_{i=1}^\infty\in X$.  Then
 $\Phi$ has a unique equilibrium state $\mu$.  The measure $\mu$ is ergodic and  has the following Gibbs property
\begin{equation}
\label{e-mat}
\mu(I)\approx \frac{\phi(I)}{\sum_{J\in \lan_n(X)} \phi(J)}\approx \exp(-nP)\phi(I),\qquad I\in \lan_n(X),\quad n\in \N.
\end{equation}
where $P=\lim_{n\to \infty}\frac{1}{n}\log \sum_{J\in \lan_n(X)} \phi(J)$. Furthermore, we have the following estimates:
 $$
\sum_{I\in \lan_n(X)}\mu(I)\log \phi(I)=n\Phi_*(\mu)+O(1),\; \sum_{I\in \lan_n(X)}\mu(I)\log \mu(I)=-nh_\mu(\sigma_X)+O(1).
$$
\end{thm}
\Proof
Define $f:\; \A^*\to [0,1]$ as in  \eqref{e-aa3}. By Lemma \ref{lem-dir},  $f\in \Omega_w(\A^*,p)$ and
$f$ satisfies $f(IJ)\preccurlyeq f(I)f(J)$ for $I,J\in \A^*$.
Let $f^*:\A^*\to [0,\infty)$ be defined as
$$f^*(I)=\sup_{n,m\geq 0}\sum_{I_1\in \A^n}\sum_{I_2\in \A^m}f(I_1II_2),\qquad I\in \A^*.
$$
Since $f(IJ)\preccurlyeq f(I)f(J)$ for $I,J\in \A^*$, we have $f^*(I)\approx f(I)$.  Hence by
Proposition \ref{pro-t1}, there exists an ergodic measure $\mu$ on $\A^\N$ such that
$\mu(I)\approx f(I)$, $I\in \A^*$.
Since $f(I)=0$ for $I\in \A^*\backslash \lan(X)$,  $\mu$ is supported on $X$.
By Lemma \ref{lem-g1}(ii), $\sum_{I\in \lan_n(X)}\phi(I)\approx \exp(nP)$, hence we have
$$
\mu(I)\approx f(I)\approx \exp(-nP)\phi(I),\quad  I\in \lan_n(X),\;n\in \N.
$$

Let $\eta$ be an ergodic equilibrium state of $\Phi$.
By Proposition \ref{pro-2.1}(i), $\Phi_*(\eta)+h_\eta(\sigma_X)=P$. By Lemma \ref{lem-ba} and \eqref{e-v2}, we have
\begin{equation}
\label{e-aa5}
\sum_{I\in \lan_n(X)}\eta(I)\log \phi(I)\geq n\Phi_*(\eta)+O(1),\; -\sum_{I\in \lan_n(X)}\eta(I)\log \eta(I)\geq nh_\eta(\sigma_X)+O(1).
\end{equation}
Thus we have
\begin{equation}
\label{e-3.10}
\begin{split}
O(1)\leq &\sum_{I\in \lan_n(X)}\Big(\eta(I)\log \phi(I)-\eta(I)\log
\eta(I)\Big)-nP\\
=&\sum_{I\in \lan_n(X)}\Big(\eta(I)\log \mu(I)-\eta(I)\log
\eta(I)\Big)+O(1).
\end{split}
\end{equation}
That is,  $\sum_{I\in \lan_n(X)}\eta(I)\log \mu(I)-\eta(I)\log
\eta(I)\geq O(1)$. By Lemma \ref{lem-4.1}, $\eta$ is absolutely continuous with respect to $\mu$. Since both $\mu$ and $\eta$ are ergodic, we have $\eta=\mu$ (cf. \cite[Theorem 6.10(iv)]{Wal82}). This implies that $\mu$ is the unique ergodic equilibrium state of $\Phi$. By Proposition \ref{pro-2.1}(ii), $\mu$ is the unique equilibrium state of $\Phi$.

Since $\eta=\mu$, by \eqref{e-3.10}, we have
\begin{equation*}
\begin{split}
\sum_{I\in \lan_n(X)}&\Big(\eta(I)\log \phi(I)-\eta(I)\log
\eta(I)\Big)-n\Phi_*(\eta)-nh_\eta(\sigma_X)\\
&=\sum_{I\in \lan_n(X)}\Big(\eta(I)\log \phi(I)-\eta(I)\log
\eta(I)\Big)-nP=O(1).
\end{split}
\end{equation*}
This together with \eqref{e-aa5} yields the estimates:
$$
\sum_{I\in \lan_n(X)}\eta(I)\log \phi(I)= n\Phi_*(\eta)+O(1),\; -\sum_{I\in \lan_n(X)}\eta(I)\log \eta(I)= nh_\eta(\sigma_X)+O(1).
$$
This completes the proof of Theorem \ref{thm-3.1}.
\qed

\begin{rem}{\rm The introduction of $\D_w(X,p)$ and $\D(X,p)$ was inspired by the work \cite{FeLa02}.
 Indeed, Theorem \ref{thm-3.1} was first setup in \cite{FeLa02} for a class of $\phi\in \D_w(X,p)$, where
$X$ is an irreducible subshift of finite type and, $\phi$ is given  by the norm of products of non-negative matrices satisfying an irreducibility condition (see \cite[Theorem 3.2]{FeLa02}, \cite[Theorem 3.1]{Fen04}). Although the approach in \cite{FeLa02} can be adapted to prove \eqref{e-mat} under our general settings,
we like to provide the above  short proof using Proposition \ref{pro-t1}. Independently, Theorem \ref{thm-3.1} was set up in \cite{Yay09a} in the special case that
$X$ is a  mixing subshift of finite type,  and $\phi$ a certain  element in $\D(X,p)$, through an approach similar to \cite{FeLa02}.
}

\end{rem}

In the end of this section, we give the following easy-checked, but important fact.
\begin{lem}
\label{lem-3.5}
Let $(X,\sigma_X)$, $(Y,\sigma_Y)$ be one-sided subshifts over  finite alphabets $\A, \A'$, respectively.
Assume that $Y$ is a factor of $X$ with a one-block factor map $\pi:\; X\to Y$.  Let $p\in \N$ and $a>0$. For $\phi\in \D_w(X,p)$, define $\phi^a:\; \lan(X)\to [0,\infty)$ and $\psi:\; \lan(Y)\to [0,\infty)$  by
$$
\phi^a(I)=\phi(I)^a \mbox{ for }I\in \lan(X),\quad \psi(J)=\sum_{I\in \lan(X):\; \pi I=J}\phi(I) \mbox{ for }J\in \lan(Y).
$$
Then $\phi^a\in \D_w(X,p)$ and $\psi\in \D_w(Y, p)$. Furthermore if  $\phi\in \D(X,p)$, then $\phi^a\in \D(X,p)$ and $\psi\in \D(Y, p)$.
\end{lem}
\Proof Clearly $\phi^a\in \D_w(X,p)$. Here we show  $\psi\in \D_w(Y, p)$.  Observe that for $J_1J_2\in \lan(Y)$,
\begin{equation*}
\begin{split}
\psi(J_1J_2)&=\sum_{I_1I_2\in \lan(X):\; \pi I_1=J_1,\;\pi I_2=J_2}\phi(I_1I_2)\\
&\leq \sum_{I_1I_2\in \lan(X):\; \pi I_1=J_1,\;\pi I_2=J_2}c^{-1}\phi(I_1)\phi(I_2)\\
&\leq \sum_{I_1\in \lan(X):\; \pi I_1=J_1}\;\sum_{I_2\in \lan(X):\;\pi I_2=J_2}c^{-1}\phi(I_1)\phi(I_2)=c^{-1}\psi(J_1)\psi(J_2).
\end{split}
\end{equation*}
Furthermore for any $J_1, J_2\in \lan(Y)$,
\begin{equation*}
\begin{split}
&\sum_{W\in \bigcup_{i=0}^p\lan_i(Y):\; J_1WJ_2\in \lan(Y)}\psi(J_1WJ_2)\\
 &\mbox{}\quad\quad =\sum_{I_1\in \lan(X):\; \pi I_1=J_1}
\sum_{I_2\in \lan(X):\; \pi I_2=J_2}\sum_{K\in \bigcup_{i=0}^p\lan_i(X):\; I_1KI_2\in \lan(X)}\phi(I_1KI_2)\\
&\mbox{}\quad\quad \geq \sum_{I_1\in \lan(X):\; \pi I_1=J_1}\;\sum_{I_2\in \lan(X):\;\pi I_2=J_2}c\phi(I_1)\phi(I_2)=c\psi(J_1)\psi(J_2).
\end{split}
\end{equation*}
Therefore there exists $W\in \bigcup_{i=0}^p\lan_i(Y)$, such that $J_1WJ_2\in \lan(Y)$, and  $\psi(J_1WJ_2)\geq \frac{c}{L}\psi(J_1)\psi(J_2)$,
 where $L$ denotes the cardinality of $\bigcup_{i=0}^p\lan_i(Y)$.  Hence  $\psi\in \D_w(Y, p)$.
A similar argument shows that $\psi\in \D(Y, p)$ whenever $\phi\in \D(X,p)$.
\qed

\section{Uniqueness of weighted equilibrium states: $k=2$}
\label{S-6}
Assume that  $(X,\sigma_X)$ is a one-sided subshift over a
finite alphabet $\A$. Let $(Y,\sigma_Y)$ be a one-sided subshift factor of $X$
with a one-block factor map $\pi:\; X\to Y$.

Let $\ba=(a_1, a_2)\in \R^2$ so that  $a_1> 0$ and $a_2\geq
0$.  Assume that $\D_w(X, p)\neq \emptyset$ for some $p\in \N$, equivalently, $X$ satisfies weak $p$-specification. Let $\phi\in \D_w(X, p)$.
Define $\phi^{(2)}:\; \lan(Y)\to [0,\infty)$  by
\begin{equation}
\label{e-6.1}
\phi^{(2)}(J)=\Big(\sum_{I\in \lan_n(X):\; \pi I=J}\phi(I)^{\frac{1}{a_1}}\Big)^{a_1} \mbox{ for } J\in \lan_n(Y),\; n\in \N.
\end{equation}
Furthermore, define $\phi^{(3)}:\; \N\to [0,\infty)$ by
\begin{equation}
\label{e-6.2}
\phi^{(3)}(n)=\sum_{J\in \lan_n(Y)}\phi^{(2)}(J)^{\frac{1}{a_1+a_2}},\quad n\in \N.
\end{equation}

The main result of this section is the following.

\begin{thm}
\label{thm-6.1} Let $\phi\in \D_w(X,p)$. Let $\Phi=(\log \phi_n)_{n=1}^\infty\in \C_{sa}(X,\sigma_{X})$ be generated by  by $\phi$, i.e. $\phi_n(x)=\phi(x_1\cdots x_n)$
for $x=(x_i)_{i=1}^\infty\in X$. Then $\Phi$ has a unique $\ba$-weighted equilibrium state $\mu$. Furthermore,
$\mu$ is ergodic and has the following properties:
\begin{itemize}
\item[{\rm (i)}] $\mu(I)\approx \widetilde{\phi}^*(I) \succcurlyeq \widetilde{\phi}(I)$ for $I\in \lan(X)$, where $\widetilde{\phi},\widetilde{\phi}^*:\;\lan(X)\to [0,\infty)$ are defined by
\begin{equation}
\label{e-6.2'}
\widetilde{\phi}(I)=
\frac{
\phi(I)^{\frac{1}{a_1}}}
{\phi^{(2)}(\pi I)^{\frac{1}{a_1}}}\cdot\frac{
\phi^{(2)}(\pi I)^{\frac{1}{a_1+a_2}}}
{\phi^{(3)}(n)},\quad I\in \lan_n(X),\; n\in \N
\end{equation}
and $$\widetilde{\phi}^*(I)=\sup_{m,n\geq 0}\sum_{I_1\in \lan_m(X),\; I_2\in \lan_n(X):\; I_1II_2\in \lan(X)}\widetilde{\phi}(I_1II_2),\quad I\in \lan(X).$$

\item[{\rm (ii)}] $\liminf_{n\to \infty} \sum_{i=0}^p\mu\left(A\cap \sigma_{X}^{-n-i}(B)\right)\succcurlyeq \mu(A)\mu(B)$ for Borel sets $A,B\subseteq X$.
\item[{\rm (iii)}]We have the estimates:
\begin{equation*}
\begin{split}
&\sum_{I\in \lan_n(X)}\mu(I)\log \mu(I)=\sum_{I\in \lan_n(X)}\mu(I)\log \widetilde{\phi}(I)+O(1)=-nh_\mu(\sigma_{X})+O(1),\\
&\sum_{I\in \lan_n(X)}\mu(I)\log \phi(I)=n\Phi_*(\mu)+O(1).
\end{split}
\end{equation*}
\end{itemize}
Moreover, if $\phi\in \D(X,p)$, then instead of (ii) we have
\begin{itemize}
\item[{\rm (iv)}]
$\liminf_{n\to \infty} \mu\left(A\cap \sigma_{X}^{-n}(B)\right)\succcurlyeq \mu(A)\mu(B)$ for Borel sets $A,B\subseteq X$.
\end{itemize}
\end{thm}

\Proof By \eqref{e-6.2'}, we have
$$
\widetilde{\phi}(I)=\frac{\phi(I)^{\frac{1}{a_1}}}{\theta(I)},\quad I\in \lan_n(X),\; n\in \N,
$$
where  $\theta(I)$ is given by
$$
\theta(I)=\phi^{(3)}(n)\phi^{(2)}(\pi I)^
{\frac{1}{a_1}-\frac{1}{a_1+a_2}},\quad I\in \lan_n(X),\; n\in \N.
$$
We claim that $\widetilde{\phi}$ and  $\theta$ satisfy the following  properties:
\begin{itemize}
\item[(a)] $\sum_{I\in \lan_n(X)}\widetilde{\phi}(I)=1$ for each $n\in \N$.
\item[(b)]For any $I\in \lan(X)$, if $\phi(I)>0$ then $\theta(I)>0$.
\item[(c)]$\theta(I_1I_2)\preccurlyeq \theta(I_1)\theta(I_2)$ for $I_1I_2\in \lan(X)$.
\end{itemize}
Property (a) follows immediately from the definition of $\widetilde{\phi}$. To see (b), one observes that if  $\phi(I)>0$ for some $I\in \lan_n(X)$,  then so are $\phi^{(2)}(\pi I)$  and  $\phi^{(3)}(n)$, hence $\theta(I)>0$. To see (c),  by Lemma \ref{lem-3.5}, $\phi^{(2)}\in \D_w(Y,p)$ and thus
$$
\phi^{(2)}(\pi (I_1I_2))\preccurlyeq \phi^{(2)}(\pi I_1)\phi^{(2)}(\pi I_2),\qquad I_1I_2\in \lan(X).\;
$$
Furthermore by Lemma \ref{lem-g1}, $\phi^{(3)}(n+m)\approx \phi^{(3)}(n)\phi^{(3)}(m)$. Hence (c) follows.

Extend $\widetilde{\phi},\; \widetilde{\phi}^*:\; \A^*\to [0,\infty)$ by setting
$\widetilde{\phi}(I)=\widetilde{\phi}^*(I)=0$ for $I\in \A^*\backslash \lan(X)$.  By (a), (b), (c) and Lemma \ref{lem-g1}(i), we see that $\widetilde{\phi}\in
\Omega_w(\A^*,p)$.  Hence by Proposition \ref{pro-t1}, there exists an ergodic measure $\mu\in \M(\A^\N,\sigma)$ such that
\begin{equation}
\label{e-6.4'}
\mu(I)\approx \widetilde{\phi}^*(I)\succcurlyeq \widetilde{\phi}(I),\quad I\in \A^\N.
\end{equation}
Moreover, $\mu$ satisfies
\begin{equation}
\label{e-6.5*}
\liminf_{n\to \infty} \sum_{i=0}^p\mu\left(A\cap \sigma^{-n-i}(B)\right)\succcurlyeq \mu(A)\mu(B) \mbox{ for Borel sets $A,B\subseteq \A^\N$}.
\end{equation}
By \eqref{e-6.4'}, $\mu$ is supported on $X$ and $\mu\in \M(X,\sigma_{X})$.

Let $\Phi^{(2)}=(\log \phi^{(2)}_n)_{n=1}^\infty\in \C_{sa}(Y,\sigma_{Y})$ be generated by $\phi^{(2)}$, i.e. $$\phi^{(2)}_n(y)=\phi^{(2)}(y_1\cdots y_n)\mbox{  for }y=(y_i)_{i=1}^\infty\in Y.$$
 Define $\widetilde{\psi}:\; \lan(Y)\to [0,\infty)$ by
$$
\widetilde{\psi}(J)=
\frac{
\phi^{(2)}(J)^{\frac{1}{a_1+a_2}}}
{\phi^{(3)}(n)},\quad J\in \lan_n(Y),\; n\in \N.
$$
 By the definitions of $\widetilde{\phi}$ and $\widetilde{\psi}$, we have
\begin{equation}
\label{e-6.3}
\widetilde{\phi}(I)=\frac{\phi(I)^{\frac{1}{a_1}}}
{\phi^{(2)}(\pi I)^{\frac{1}{a_1}}}\cdot \widetilde{\psi}(\pi I), \quad I\in \lan(X).
\end{equation}

Since $\phi^{(2)}\in \D_w(Y,p)$, by Lemma \ref{lem-3.5},  $(\phi^{(2)})^{1/(a_1+a_2)}\in \D_w(Y,p)$. Hence by Theorem \ref{thm-3.1}, $\frac{1}{a_1+a_2}\Phi^{(2)}$ has a unique equilibrium state $\nu\in \M(Y,\sigma_{Y})$ and $\nu$ satisfies the properties
\begin{equation}
\label{e-6.4}
\sum_{J\in \lan_n(Y)}\nu(J)\log \nu(J)=\sum_{J\in \lan_n(Y)}\nu(J)\log \widetilde{\psi}(J)+O(1)=-nh_\nu(\sigma_{Y})+O(1),
\end{equation} and
\begin{equation}
\label{e-6.4''}
\sum_{J\in \lan_n(Y)}\nu(J)\log \phi^{(2)}(J)=n\Phi^{(2)}_*(\nu)+O(1).
\end{equation}

Assume that $\eta$ is an ergodic $\ba$-equilibrium state of $\Phi$. By Corollary \ref{cor-1.1}(i), $\eta\circ \pi^{-1}=\nu$ and
$\eta$ is a conditional equilibrium state of $\frac{1}{a_1}\Phi$ with respect to $\nu$, that is,
\begin{equation}
\label{e-6.5}
\frac{1}{a_1}\Phi_*(\eta)+h_\eta(\sigma_{X})-h_\nu(\sigma_{Y})=\frac{1}{a_1}\Phi_*^{(2)}(\nu).
\end{equation}
By \eqref{e-6.4} and \eqref{e-6.4''},  we have
\begin{equation}
\label{e-6.5'}
nh_\nu(\sigma_{Y}) + \frac{n}{a_1} \Phi_*^{(2)} (\nu)=-\sum_{J\in \lan_n(Y)} \nu(J) \log \frac{\widetilde{\psi}(J)} {\phi^{(2)}(J)^{\frac{1}{a_1}}}+O(1).
\end{equation}
By Lemma \ref{lem-ba}(ii) and \eqref{e-v2}, we have
\begin{equation}
\label{e-6.5''}
\sum_{I\in \lan_n(X)}\eta(I)\log \phi(I)\geq n\Phi_*(\eta)+O(1),\quad
-\sum_{I\in \lan_n(X)} \eta(I)\log \eta(I)\geq nh_\eta(\sigma_{X}).
\end{equation}
Combining \eqref{e-6.5}, \eqref{e-6.5'} and \eqref{e-6.5''}, we obtain
\begin{equation}
\label{e-6.10}
\begin{split}
O(1)\leq &\sum_{I\in \lan_n(X)}\Big(\eta(I)\log (\phi(I)^{\frac{1}{a_1}})-\eta(I)\log
\eta(I)\Big)-\frac{n}{a_1}\Phi_*(\eta)-nh_\eta(\sigma_X)\\
=&\sum_{I\in \lan_n(X)}\Big(\eta(I)\log (\phi(I)^{\frac{1}{a_1}})-\eta(I)\log
\eta(I)\Big)-nh_\nu(\sigma_{X_2})-\frac{n}{a_1}\Phi_*^{(2)}(\nu)\\
=& \sum_{I\in \lan_n(X)}\Big(\eta(I)\log (\phi(I)^{\frac{1}{a_1}})-\eta(I)\log
\eta(I)\Big)\\
&\qquad\qquad\qquad +\sum_{J\in \lan_n(Y)} \nu(J) \log \frac{\widetilde{\psi}(J)} {\phi^{(2)}(J)^{\frac{1}{a_1}}}+O(1)\\
=&\sum_{I\in \lan_n(X)}\Big(\eta(I)\log \frac
 {\phi(I)^{\frac{1}{a_1}}\widetilde{\psi}(\pi I)}{\phi^{(2)}(\pi I)^{\frac{1}{a_1}}}
 -\eta(I)\log
\eta(I)\Big)
 + O(1)\\
 =& \sum_{I\in \lan_n(X)}\Big(\eta(I)\log \widetilde{\phi}(I)-\eta(I)\log
\eta(I)\Big)+ O(1)\qquad (\mbox{by \eqref{e-6.3}}).
\end{split}
\end{equation}
That is,
\begin{equation}
\label{e-6.11}
 \sum_{I\in \lan_n(X)}\Big(\eta(I)\log \widetilde{\phi}(I)-\eta(I)\log
\eta(I)\Big)\geq O(1).
\end{equation}
Combining \eqref{e-6.11} and \eqref{e-6.4'} yields
\begin{equation}
\label{e-6.11'}
\begin{split}
 \sum_{I\in \lan_n(X)}&\Big(\eta(I)\log \mu(I)-\eta(I)\log
\eta(I)\Big)\\
&\geq \sum_{I\in \lan_n(X)}\Big(\eta(I)\log \widetilde{\phi}(I)-\eta(I)\log
\eta(I)\Big)+O(1)\geq O(1).
\end{split}
\end{equation}
By \eqref{e-6.11'} and Lemma \ref{lem-4.1}, $\eta$ is absolutely continuous with respect to $\mu$. Since both $\mu$ and $\eta$ are ergodic, we have $\eta=\mu$ (cf. \cite[Theorem 6.10(iv)]{Wal82}). This implies that $\mu$ is the unique ergodic $\ba$-weighted equilibrium state of $\Phi$. By Corollary \ref{cor-1.1}(iii), $\mu$ is the unique  $\ba$-weighted equilibrium state of $\Phi$. Now parts (i), (ii) of the theorem follow from \eqref{e-6.4'}-\eqref{e-6.5*}.

To show (iii),  due to  $\eta=\mu$, the left hand side of \eqref{e-6.11'} equals $0$. Hence by \eqref{e-6.11'},
\begin{equation}
\label{e-6.15}
\sum_{I\in \lan_n(X)}\Big(\eta(I)\log \widetilde{\phi}(I)-\eta(I)\log
\eta(I)\Big)=O(1).
\end{equation}
Combining \eqref{e-6.15} and \eqref{e-6.10} yields
\begin{equation}
\label{e-6.16}
\sum_{I\in \lan_n(X)}\Big(\eta(I)\log (\phi(I)^{\frac{1}{a_1}})-\eta(I)\log
\eta(I)\Big)-\frac{n}{a_1}\Phi_*(\eta)-nh_\eta(\sigma_X)=O(1).
\end{equation}
 However \eqref{e-6.16} and \eqref{e-6.5''} imply
\begin{equation}
\label{e-6.17}
\sum_{I\in \lan_n(X)}\eta(I)\log \phi(I)= n\Phi_*(\eta)+O(1),\;
-\sum_{I\in \lan_n(X)} \eta(I)\log \eta(I)= nh_\eta(\sigma_{X})+O(1).
\end{equation}
Now part (iii) follows from \eqref{e-6.17} and \eqref{e-6.15}. To see (iv), note that whenever $\phi\in \D(X,p)$,
we have $\widetilde{\phi}\in \Omega(\A^*, p)$, following from (a)-(c). Now (iv) follows from Proposition \ref{pro-t1}(vi).
  This finishes the proof of the theorem.
\qed

\section{Uniqueness of weighted equilibrium states: $k\geq 2$}
\label{S-7}
Let $k\geq 2$ be an integer. Assume that $(X_i, \sigma_{X_i})$ ($i=1,\ldots,k)$ are one-sided subshifts over finite alphabets so that
$X_{i+1}$ is a factor of $X_i$ with a one-block factor map $\pi_i:\; X_i\to X_{i+1}$ for $i=1,\ldots,k-1$.
For convenience, we use $\pi_0$ to denote the identity map on $X_1$. Define
 $\tau_i:\;X_1\to X_{i+1}$ by $\tau_i=\pi_i\circ\pi_{i-1}\circ
\cdots \circ \pi_0$ for $i=0,1,\ldots,k-1$.

Let $\ba=(a_1,\ldots, a_k)\in \R^k$ so that  $a_1> 0$ and $a_{i}\geq
0$ for $i>1$.  Let $\phi\in \D_w(X_1, p)$. Set $\phi^{(1)}=\phi$ and define $\phi^{(i)}:\; \lan(X_i)\to [0,\infty)$ ($i=2,\ldots,k$) recursively by
\begin{equation*}
\label{e-7.1}
\phi^{(i)}(J)=\Big(\sum_{I\in \lan_n(X_{i-1}):\; \pi_{i-1} I=J}\phi^{(i-1)}(I)^{\frac{1}{a_1+\cdots+a_{i-1}}}\Big)^{a_1+\cdots+a_{i-1}}
\end{equation*}
for $n\in \N$, $J\in \lan_n(X_{i})$.
Furthermore, define $\phi^{(k+1)}:\; \N\to [0,\infty)$ by
\begin{equation*}
\label{e-7.2}
\phi^{(k+1)}(n)=\sum_{I\in \lan_n(X_{k})}\phi^{(k)}(I)^{\frac{1}{a_1+\cdots+a_k}}.
\end{equation*}

Let $\Phi=(\log \phi_n)_{n=1}^\infty\in \C_{sa}(X_1,\sigma_{X_1})$ be generated by $\phi$. Say that $\mu\in \M(X_1,\sigma_{X_1})$ is an {\it $\ba$-weighted equilibrium state of $\Phi$} if
$$
\Phi_*(\mu)+\sum_{i=1}^k a_ih_{\mu\circ \tau_{i-1}^{-1}}(\sigma_{X_i})=\sup_{\eta\in \M(X,\sigma_{X_1})}\left(\Phi_*(\eta)+\sum_{i=1}^k a_ih_{\eta\circ \tau_{i-1}^{-1}}(\sigma_{X_i})\right).
$$
Let $\I(\Phi,\ba)$ be the collection of all {\it $\ba$-weighted equilibrium state of $\Phi$}.

Let $\Phi^{(2)}\in \C_{sa}(X_2,\sigma_{X_2})$ be generated by $\phi^{(2)}$. By a proof essentially identical  to that of Corollary \ref{cor-1.1},
we have
\begin{lem}
\label{lem-t}
\begin{itemize}
\item[{\rm (i)}]
$\I(\Phi,\ba)$ is a non-empty compact convex subset of $\M(X_1,\sigma_{X_1})$. Each extreme point of $\I(\Phi,\ba)$ is ergodic.
\item[{\rm (ii)}] $\mu\in \I(\Phi,\ba)$ if and only if $\mu\in \I_{\mu\circ \pi_1^{-1}}(\frac{1}{a_1}\Phi)$ together with $\mu\circ \pi^{-1}\in \I(\Phi^{(2)}, \bb)$, where $\bb=(a_1+a_2, a_3,\ldots, a_k)\in \R^{k-1}$.
\item[{\rm (iii)}]  $\I(\Phi,\ba)$ is a singleton if and only if $\I(\Phi^{(2)}, \bb)$ is a singleton $\{\nu\}$ and, $\I_{\nu}(\frac{1}{a_1}\Phi)$ contains a unique ergodic measure.
\end{itemize}
\end{lem}

As the high dimensional version of Theorem \ref{thm-6.1}, we have

\begin{thm}
\label{thm-7.1}
Let $\phi\in \D_w(X_1, p)$. Let $\Phi=(\log \phi_n)_{n=1}^\infty\in \C_{sa}(X_1,\sigma_{X_1})$ be generated by $\phi$. Then $\Phi$ has a unique $\ba$-weighted equilibrium state $\mu$. Furthermore,
$\mu$ is ergodic and has the following properties:
\begin{itemize}
\item[{\rm (i)}] $\mu(I)\approx \widetilde{\phi}^*(I)\succcurlyeq \widetilde{\phi}(I)$ for $I\in \lan(X_1)$, where $\widetilde{\phi},\; \widetilde{\phi}^*:\; \lan(X_1)\to [0,\infty)$ are defined respectively by
\begin{equation}
\label{e-7.2'}
\widetilde{\phi}(I)=
\left(\prod_{i=1}^{k-1}
\frac{
\phi^{(i)}(\tau_{i-1}I)^{\frac{1}{a_1+\cdots+a_i}}}
{\phi^{(i+1)}(\tau_{i}I)^{\frac{1}{a_1+\cdots+a_i}}}\right)\cdot\frac{
\phi^{(k)}(\tau_{k-1}I)^{\frac{1}{a_1+\cdots+a_k}}}
{\phi^{(k+1)}(n)}
\end{equation}
for $I\in \lan_n(X_1)$, $n\in \N$, and
$$\widetilde{\phi}^*(I)=\sup_{m,n\geq 0}\sum_{I_1\in \lan_m(X_1),\; I_2\in \lan_n(X_1):\; I_1II_2\in \lan(X_1)}\widetilde{\phi}(I_1II_2),\quad I\in \lan(X_1).$$
\item[{\rm (ii)}] $\liminf_{n\to \infty} \sum_{i=0}^p\mu\left(A\cap \sigma_{X_1}^{-n-i}(B)\right)\succcurlyeq \mu(A)\mu(B)$ for Borel sets $A,B\subseteq X_1$.
\item[{\rm (iii)}]We have the estimates:
\begin{equation*}
\begin{split}
&\sum_{I\in \lan_n(X_1)}\mu(I)\log \mu(I)=\sum_{I\in \lan_n(X_1)}\mu(I)\log \widetilde{\phi}(I)+O(1)=-nh_\mu(\sigma_{X_1})+O(1),\\
&\sum_{I\in \lan_n(X_1)}\mu(I)\log \phi(I)=n\Phi_*(\mu)+O(1).
\end{split}
\end{equation*}
\end{itemize}
Moreover, if $\phi\in \D(X_1,p)$, then instead of (ii) we have
\begin{itemize}
\item[{\rm (iv)}]
$\liminf_{n\to \infty} \mu\left(A\cap \sigma_{X_1}^{-n}(B)\right)\succcurlyeq \mu(A)\mu(B)$ for Borel sets $A,B\subseteq X_1$.
\end{itemize}
\end{thm}
\Proof We prove the theorem by induction on the dimension $k$.
By Theorem \ref{thm-6.1}, Theorem \ref{thm-7.1} is true when the dimension equals $2$.
Now assume that the theorem is true when the dimension equals  $k-1$.  In the following we prove that the theorem is also true when the dimension equals  $k$.

By \eqref{e-7.2'}, we have
$$
\widetilde{\phi}(I)=\frac{\phi(I)^{\frac{1}{a_1}}}{\theta(I)},\quad I\in \lan_n(X_1),\; n\in \N,
$$
where  $\theta(I)$ is given by
$$
\theta(I)=\phi^{(k+1)}(n)\prod_{i=2}^{k}\phi^{(i)}(\tau_{i-1}I)^
{\frac{1}{a_1+\cdots+a_{i-1}}-\frac{1}{a_1+\cdots+a_{i}}},\quad I\in \lan_n(X_1),\; n\in \N.
$$
By Lemma \ref{lem-3.5} and Lemma \ref{lem-g1}(ii), we have $\phi^{(i)}\in \D_w(X_i,p)$ for $i=2,\ldots,k$, and
$\phi^{(k+1)}(n+m)\approx \phi^{(k+1)}(n)\phi^{(k+1)}(m)$. Similar to the proof of Theorem \ref{thm-6.1}, we can show that  $\widetilde{\phi}$ and  $\theta$ satisfy the following  properties:
\begin{itemize}
\item[(a)] $\sum_{I\in \lan_n(X_1)}\widetilde{\phi}(I)=1$ for each $n\in \N$.
\item[(b)]For any $I\in \lan(X_1)$, if $\phi(I)>0$ then $\theta(I)>0$.
\item[(c)]$\theta(I_1I_2)\preccurlyeq \theta(I_1)\theta(I_2)$ for $I_1I_2\in \lan(X_1)$.
\end{itemize}
Extend $\widetilde{\phi},\; \widetilde{\phi}^*:\; \A^*_1\to [0,\infty)$ by setting
$\widetilde{\phi}(I)=\widetilde{\phi}^*(I)=0$ for $I\in \A^*_1\backslash \lan(X_1)$.  By (a), (b), (c) and Lemma \ref{lem-g1}(i), we see that $\widetilde{\phi}\in
\Omega_w(\A_1^*,p)$.  Hence by Proposition \ref{pro-t1}, there exists an ergodic measure $\mu\in \M(\A_1^\N,\sigma)$ such that
\begin{equation}
\label{e-7.4'}
\mu(I)\approx \widetilde{\phi}^*(I)\succcurlyeq \widetilde{\phi}(I),\quad I\in \A^\N_1.
\end{equation}
Moreover, $\mu$ satisfies
\begin{equation*}
\label{e-7.5*}
\liminf_{n\to \infty} \sum_{i=0}^p\mu\left(A\cap \sigma^{-n-i}(B)\right)\succcurlyeq \mu(A)\mu(B) \mbox{ for Borel sets $A,B\subseteq \A_1^\N$}.
\end{equation*}
By \eqref{e-7.4'}, $\mu$ is supported on $X_1$ and $\mu\in \M(X_1,\sigma_{X_1})$.

Let $\Phi^{(2)}=(\log \phi^{(2)}_n)_{n=1}^\infty\in \C_{sa}(X_2,\sigma_{X_2})$ be generated by $\phi^{(2)}$, i.e.  $$\phi^{(2)}_n(x)=\phi^{(2)}(x_1\cdots x_n)\mbox{  for }x=(x_i)_{i=1}^\infty\in X_2.$$
Let $\bb=(a_1+a_2, a_3,\ldots,a_k)\in \R^{k-1}$. Define $\widetilde{\psi}:\; \lan(X_2)\to [0,\infty)$ by
$$
\widetilde{\psi}(J)=
\left(\prod_{i=2}^{k-1}
\frac{
\phi^{(i)}(\xi_{i-1}J)^{\frac{1}{a_1+\cdots+a_i}}}
{\phi^{(i+1)}(\xi_{i}J)^{\frac{1}{a_1+\cdots+a_i}}}\right)\cdot\frac{
\phi^{(k)}(\xi_{k-1}J)^{\frac{1}{a_1+\cdots+a_k}}}
{\phi^{(k+1)}(n)},\quad J\in \lan(X_2),\; n\in \N,
$$
where $\xi_{1}:=Id$, and $\xi_{i}=\pi_{i}\circ \cdots\circ \pi_2$ for $i\geq 2$. By the definitions of $\widetilde{\phi}$ and $\widetilde{\psi}$, we have
\begin{equation}
\label{e-7.3}
\widetilde{\phi}(I)=\frac{\phi^{(1)}(I)^{\frac{1}{a_1}}}
{\phi^{(2)}(\pi_1I)^{\frac{1}{a_1}}}\cdot \widetilde{\psi}(\pi_1 I), \quad I\in \lan(X_1).
\end{equation}

Since $\phi^{(2)}\in \D_w(X_2,p)$, by the assumption of the induction, $\Phi^{(2)}$ has a unique $\bb$-weighted equilibrium state $\nu\in \M(X_2,\sigma_{X_2})$ and $\nu$ satisfies the properties
\begin{equation}
\label{e-7.4}
\sum_{J\in \lan_n(X_2)}\nu(J)\log \nu(J)=\sum_{J\in \lan_n(X_2)}\nu(J)\log \widetilde{\psi}(J)+O(1)=-nh_\nu(\sigma_{X_2})+O(1),
\end{equation} and
\begin{equation}
\label{e-7.4''}
\sum_{J\in \lan_n(X_2)}\nu(J)\log \phi^{(2)}(J)=n\Phi^{(2)}_*(\nu)+O(1).
\end{equation}

Assume that $\eta$ is an ergodic $\ba$-equilibrium state of $\Phi$. By Lemma \ref{lem-t}, $\eta\circ \pi^{-1}_1=\nu$ and
$\eta$ is a conditional equilibrium state of $\frac{1}{a_1}\Phi$ with respect to $\nu$, that is,
\begin{equation}
\label{e-7.5}
\frac{1}{a_1}\Phi_*(\eta)+h_\eta(\sigma_{X_1})-h_\nu(\sigma_{X_2})=\frac{1}{a_1}\Phi_*^{(2)}(\nu).
\end{equation}

Using  \eqref{e-7.4'}, \eqref{e-7.3}, \eqref{e-7.4}-\eqref{e-7.5}, and  taking a process the same as in the proof of Theorem \ref{thm-6.1},  we  prove Theorem \ref{thm-7.1} when the dimension equals $k$.
\qed

\begin{rem}
{\rm
Let $\widetilde{\phi}$ be defined as in \eqref{e-7.2'}, and let $(\eta_n)$ be a sequence in $\M(X)$ so that
$\eta_n(I)=\widetilde{\phi}(I)$ for each $I\in \lan_n(X_1)$. Then by Proposition \ref{pro-t1}(v) and the above proof,  the measure $\mu$ in Theorem \ref{thm-7.1} satisfies
$$
\mu=\lim_{n\to \infty}\frac{1}{n}\sum_{i=0}^{n-1}\eta_n\circ \sigma_{X_1}^{-i}.
$$

}
\end{rem}

\bigskip
\noindent {\it Proof of Theorem \ref{thm-1}}.  We first consider the case that $X_i$ ($i=1,\ldots,k)$ are one-sided subshifts. Recoding  $X_{k-1}$, $X_{k-1},\ldots, X_1$ recursively through their higher block representations (cf. Proposition 1.5.12 in \cite{LiMa95}), if necessary, we may assume that $\pi_i:\; X_{i}\to X_{i+1}$ ($i=1,\ldots, k-1$) are all one-block factor maps.  Recall that $X_1$ satisfies weak specification.  (Notice that this property  is preserved by recoding via  higher  block representations). Let $f\in V(\sigma_{X_1})$ (see \eqref{e-sn1} for the definition). Define
$\phi:\; \lan(X_1)\to [0,\infty)$ by
$$
\phi(I)=\sup_{x\in X_1\cap[I]}\exp(S_nf(x)),\quad I\in \lan_n(X_1),\; n\in \N,
$$
where $S_nf$ is defined as in \eqref{e-sn}.
Since $f\in V(\sigma_{X_1})$, it is direct to check that $\phi\in \D_w(X_1,p)$, where $p$ is any integer so that $X_1$ satisfies weak $p$-specification. Let $\Phi=(\log \phi_n)_{n=1}^\infty\in \C_{sa}(X_1,\sigma_{X_1})$ be generated by $\phi$.  Again by  $f\in V(\sigma_{X_1})$, we have
$\Phi_*(\mu)=\mu(f)$ for any $\mu\in \M(X_1,\sigma_{X_1})$. It follows that $\mu$ is an $\ba$-weighted equilibrium state of $f$ if and only if that,  $\mu$ is an $\ba$-weighted equilibrium state of $\Phi$. Now the theorem  follows from Theorem \ref{thm-7.1}.

Next we  consider  the case that $X_i$'s are two-sided subshifts over  finite alphabets $\A_i$'s. Again we may assume that $\pi_i$'s are one-block factor maps.  Define for $i=1,\ldots,k$,
$$
X_i^+:=\left\{(x_j)_{j=1}^\infty\in \A_i^\N:\; \exists \;(y_j)_{j\in \Z}\in X_i \mbox{ such that }x_j=y_j\mbox { for }j\geq 1\right\}.
$$
Then $(X_i^{+},\sigma_{X_i^{+}})$ becomes a one-sided subshift for each $i$. Furthermore define
$\Gamma_i:\; X_i\to X_i^{+}$ by $(x_j)_{j\in \Z}\mapsto (x_j)_{j\in \N}$. Then for each $1\leq i\leq k$,
the mapping $\mu\mapsto \mu\circ {\Gamma_i}^{-1}$ is an homeomorphism from $\M(X_i,\sigma_{X_i})$ to  $\M(X_i^+,\sigma_{X_i^+})$ which preserves the measure theoretic entropy.  Now $\pi_i:\; X_i^+\to X_{i+1}^+$ becomes
a one-block factor between one-sided subshifts for $i=1,\ldots,k-1$.    Let $f\in  V(\sigma_{X_1})$. Define
$\phi:\; \lan(X_1^+)\to [0,\infty)$ by
$$
\phi(I)=\sup_{x\in X_1:\; x_1\ldots x_n=I}\exp(S_nf(x)),\quad I\in \lan_n(X_1^+),\; n\in \N.
$$
Similarly, $\phi\in \D_w(X_1,p)$ for some $p\in \N$. Let $\Phi=(\log \phi_n)_{n=1}^\infty\in \C_{sa}(X_1^+,\sigma_{X_1^+})$ be generated by $\phi$. Due to $f\in  V(\sigma_{X_1})$, we have
$$
\mu(f)=\Phi_*(\mu\circ \Gamma_1^{-1}),\qquad \mu\in \M(X_1,\sigma_{X_1}).
$$
It follows that  $\mu$ is an $\ba$-weighted equilibrium state of $f$ if and only if that, $\mu\circ \Gamma_1^{-1}$ is an $\ba$-weighted equilibrium state of $\Phi$. Thus the results of the theorem follow from  Theorem \ref{thm-7.1}.
\qed

\noindent {\bf Acknowledgements}.
 The author is grateful to Eric Olivier for some helpful discussions about the question of Gatzouras and Peres. He was  partially
supported by the RGC grant (project CUHK401008) in the Hong Kong Special Administrative
Region, China.

\end{document}